\documentclass[twoside,12pt]{amsart}
\usepackage[T1]{fontenc}
\usepackage{amsfonts}
\usepackage{amssymb}
\usepackage{amsthm}
\usepackage{amsmath}
\usepackage{graphicx}
\usepackage{url}
\usepackage[all]{xy}
\usepackage{lscape}
\theoremstyle{plain}

\newtheorem*{theorem*}{Theorem}

\newcommand{\GU}{\mbox{\rm GU}}
\newcommand{\Sp}{\mbox{\rm Sp}}
\newcommand{\SO}{\mbox{\rm SO}}
\newcommand{\GL}{\mbox{\rm GL}}

\newcommand{\End}{\mbox{\rm End}}

\newcommand{\Infl}{\mbox{\rm Infl}}

\newcommand{\Res}{\mbox{\rm Res}}

\newcommand{\bc}{\mathbf{c}}
\newcommand{\bd}{\mathbf{d}}
\newcommand{\bs}{\mathbf{s}}

\newcommand{\cont}{\text{cont}}
\newcommand{\res}{\text{res}}
\newcommand{\wt}{\text{wt}}
\newcommand{\mand}{\quad\text{and}\quad}

\newenvironment{prf}{{\bf Proof.}}{\hfill $\Box$ \\[-1.0ex]}

\newtheorem*{define*}{Definition}

\newtheorem*{thm*}{Theorem}

\newtheorem*{lem*}{Lemma}

\newtheorem*{prp*}{Proposition}

\newtheorem*{cor*}{Corollary}

\newtheorem*{conj*}{Conjecture}

\newtheorem*{xmpl*}{Example}

\newtheorem*{rem*}{Remark}

\begin{document}
\date{\today}
\title{Harish-Chandra series in finite unitary groups and crystal graphs}
\author{Thomas Gerber, Gerhard Hiss and Nicolas Jacon}


\address{TG: Laboratoire de Math{\'e}matiques et Physique Th{\'e}orique, 
(UMR 7350, CNRS -- Universit{\'e} de Tours)
Parc de Grandmont, 37200, Tours}
\address{GH: Lehrstuhl D f\"ur Mathematik, RWTH Aachen University,
52056 Aachen, Germany}
\address{NJ: Laboratoire de Math{\'e}matique de Reims, 51687 Reims cedex 2}

\email{thomas.gerber@lmpt.univ-tours.fr}
\email{gerhard.hiss@math.rwth-aachen.de}
\email{nicolas.jacon@univ-reims.fr}

\subjclass[2000]{20C33, 20C08, 20G42, 17B37, 81R50}
\keywords{Harish-Chandra series, unitary group, branching graph, Fock space, 
crystal basis, crystal graph}

\begin{abstract}
The distribution of the unipotent modules (in non-defining prime characteristic) 
of the finite unitary groups into Ha\-rish-Chandra series is investigated. We
formulate a series of conjectures relating this distribution with the crystal
graph of an integrable module for a certain 
quantum group. Evidence for our conjectures is presented, as well as proofs for 
some of their consequences for the crystal graphs involved. In the course of our 
work we also generalize Harish-Chandra theory for some of the finite classical 
groups, and we introduce their Harish-Chandra branching graphs.
\end{abstract}

\maketitle

\section{Introduction}
\markright{HARISH-CHANDRA SERIES AND CRYSTAL GRAPHS}

Harish-Chandra theory provides a means of labelling the simple modules of a 
finite group~$G$ of Lie type in non-defining characteristics, including~$0$. 
The set of simple modules of~$G$ (up to isomorphism) is partitioned into 
disjoint subsets, the Harish-Chandra series, each arising from a cuspidal 
simple module of a Levi subgroup of~$G$. Inside each series, the modules 
are classified by the simple modules of an Iwahori-Hecke algebra arising 
from the the cuspidal module which representing the series.

This yields, however, a rather indirect labelling of the simple modules, as 
it requires the classification of the cuspidal simple modules. Moreover, for
each of these, the corresponding Iwahori-Hecke algebra has to be computed and 
its simple modules have to be classified. This program has been completed 
successfully by Lusztig for modules over fields of characteristic~$0$ 
(see~\cite{luszbuch}). For modules over fields of positive characteristic, 
only partial results are known.

In some cases a different labelling of the simple modules of~$G$ is known. 
This arises from Lusztig's classification of the simple modules in 
characteristic~$0$, together with sufficient knowledge of Brauer's theory of 
decomposition numbers. This applies in particular to the general linear groups 
$G = \GL_n(q)$ and the general unitary groups $G = \GU_n(q)$, where the 
unipotent modules (in any non-defining characteristic) are labelled by 
partitions of~$n$. For characteristic~$0$ this result is due to 
Lusztig and Srinivasan \cite{LuSri}, for prime 
characteristic it follows from work of Dipper \cite{DipTri} and Geck 
\cite{TriangularShape}. In these cases it is natural to ask how to determine
the partition of the unipotent modules into Harish-Chandra series from these 
labels of the unipotent modules, i.e.\ from the partitions of~$n$. 

By work of Dipper and Du (see \cite[Section~$4$]{DDu2}), this can be done
for the general linear groups. First attempts to find a similar description for 
the unitary groups are described in~\cite{GeHiMaI}. It turned out, however, that 
this is possible only in a favourable case, the case of linear characteristic
(see \cite[Corollary~$8.11$]{GruHi} in conjunction with the above mentioned
results by Dipper and Du). The general description of the Harish-Chandra series 
of the unitary groups and other classical groups is still open.

In this paper we present a series of conjectures which, when true and proved,
will solve generalized versions of this problem, at least for large 
characteristics. 

Let us now describe our main results and conjectures. As above,~$G$ denotes a 
finite group of Lie type, viewed as group with a split $BN$-pair of 
characteristic~$p$. We also let~$\ell$ be a prime different from~$p$. In this
introduction, by a simple module for~$G$ we will always mean an absolutely 
simple module over a field of characteristic~$0$ or~$\ell$. In Section~$2$ we 
introduce a generalization of Harish-Chandra theory if~$G$ is a unitary, 
symplectic and odd dimensional orthogonal groups. Thus the Weyl group of~$G$,
as group with a $BN$-pair, is of type~$B$.
Instead of using all Levi subgroups for Harish-Chandra induction, we restrict 
to what we call pure Levi subgroups: those that arise from a connected subset
of the Dynkin diagram of~$G$ which is either empty or else contains the first
node ajacent to the double edge.
This way we obtain more cuspidal modules, which we call weakly cuspidal.
All main results of Harish-Chandra theory remain valid in this more 
general context. In particular, we obtain a distribution of the simple
modules into weak Harish-Chandra series (Proposition~\ref{WeakHCSeries}).
The usual Harish-Chandra series are unions of weak Harish-Chandra series.
In characteristic~$0$, the two notions coincide for unipotent modules, as
a Levi subgroup having a unipotent cuspidal module is pure by Lusztig's
classification.

In Section~$3$ we prove some results on the endomorphism ring of a 
Harish-Chandra induced weakly cuspidal module. Theorem~\ref{StructureHLX} 
states that, under some mild restrictions, this endomorphism ring is in fact an 
Iwahori-Hecke algebra of type~$B$. Some information about the parameters of this
algebra are also given. For example, if a simple weakly cuspidal module in 
characteristic~$\ell$ lies in a block containing an ordinary cuspidal module, 
then the parameters of the two Iwahori-Hecke algebras are related through 
reduction modulo~$\ell$.

In Section~$4$ we define the Harish-Chandra branching graph for the
unipotent modules of the classical groups considered. This graph re\-cords
the socle composition factors of Harish-Chandra induced unipotent modules, 
very much in the spirit of Kleshchev's branching rules for modules of 
symmetric groups (see \cite{KleshchevI,KleshchevII,KleshchevIII, KleshchevIV},
in particular \cite[Theorem~$0.5$]{KleshchevII}).

Section~$5$ contains our conjectures. These are restricted to the case of
the unitary groups. We thus let $G = \GU_n(q)$ from now on and we write~$e$ 
for the multiplicative order of~$-q$ in a field of characteristic~$\ell$. 
Following \cite[Definition~$5.3$]{GruHi}, we call~$\ell$ linear for~$G$, 
if~$e$ is even. For our conjectures, however,  we assume that~$e$ is odd
and larger than~$1$, so that in particular~$\ell$ is non-linear for~$G$. 
(The case $e = 1$, i.e.\ $\ell \mid q + 1$ has been settled in~\cite{GeHiMaII}.)
Conjecture~\ref{Compatibility} 
concerns the relation between Harish-Chandra series of ordinary modules and 
those in characteristic~$\ell$. It predicts that if two unipotent modules 
of~$G$, labelled by the partitions $\lambda$ and $\mu$, respectively, lie in 
the same weak Harich-Chandra series, then~$\lambda$ and~$\mu$ have the same 
$2$-core, i.e.\ the ordinary unipotent modules labelled by these two 
partitions also lie in the same Harish-Chandra series. 
In this sense the $\ell$-modular Harish-Chandra series (of unipotent modules)
form a refinement of the ordinary Harish-Chandra series. According to
Conjecture~\ref{CuspidalityAndECore}, the $e$-core of~$\lambda$ should be a 
$2$-core, if~$\lambda$ labels a weakly cuspidal unipotent module. This
amounts to the assertion that if a unipotent $\ell$-block contains a weakly
cuspidal module, then the block also contains an ordinary cuspidal module
(not necessarily unipotent). Conjecture~\ref{GraphComparison} relates the 
Harish-Chandra branching graphs with crystal graphs arising from canonical 
bases in submodules of Fock spaces of level~$2$, which are acted on by the 
quantum group $\mathcal{U}'_v( \widehat{\mathfrak{sl}_e} )$. This is in analogy 
to the case of Kleshchev's branching graph in characteristic~$p$, which is 
isomorphic to the crystal graph of a Fock space of level~$1$ with an action of 
the quantum group $\mathcal{U}'_v( \widehat{\mathfrak{sl}_p} )$ (see
\cite{KleshchevI,KleshchevII,KleshchevIII, KleshchevIV}). The conjecure is also 
put in perspective by the results of Shan \cite{Shan} on the branching rules 
on the category $\mathcal{O}$ of the cyclotomic rational double affine Hecke 
algebras.
Finally, 
Conjecture~\ref{Uglov} 
is just a weaker form of Conjecture~\ref{GraphComparison}. Its statement gives 
an algorithm to compute the distribution of the unipotent modules in 
characteristic~$\ell$ into weak Harish-Chandra series from the combinatorics of 
the crystal graph involved. In our conjectures we assume that $\ell$ is large 
enough (compared to~$n$), without specifying any bound. In the computed 
examples, $\ell > n$ is good enough.

In Section~$6$ we collect our evidence for the conjectures. 
In Theorem~\ref{MainSeriesLabelling} we prove that Conjecture~\ref{Uglov} holds 
for some subgraphs of the Harish-Chandra branching graph and the crystal graph, 
respectively. It is a generalization of the main result of Geck~\cite{GeckPS} 
for principal series to other ordinary Harish-Chandra series. Similarly, 
Theorem~\ref{WeightOne} asserts that parts of our conjectures hold for
blocks of weight~$1$, i.e.\ blocks with cyclic defect groups.
We also compute the parameters of the Iwahori-Hecke algebra corresponding to 
a weakly cuspidal module under the assumption that 
Conjecture~\ref{CuspidalityAndECore} holds true
(Proposition~\ref{ParameterOf2Core}). Finally, the truth of 
Conjectures~\ref{GraphComparison} and~\ref{Uglov} implies an isomorphism
of certain connected components of crystal graphs with different
parameters. This is discussed in~\ref{DifferentGraphs}.

In Section~$7$ we prove that the consequences implied by the conjectures for 
the crystal graphs are indeed true. This adds more evidence to our conjectures. 
Conjecture~\ref{Uglov} implies that a weakly cuspidal module is labelled by a 
partition which gives rise to a highest weight vertex in the crystal graph. Such 
partitions can be characterized combinatorially (see \cite{JL}). We prove in 
Theorem~\ref{ConsequenceConjectureUglov} that the corresponding $e$-core is 
indeed a $2$-core, as predicted by Conjecture~\ref{CuspidalityAndECore}. 
In \cite[Theorem~$8.3$]{GeHiMaI} we had proved that the unipotent module of~$G$
labelled by the partition $(1^n)$ is cuspidal if and only if~$\ell$ divides~$n$ 
or~$n-1$. We prove that the anologous statement holds for corresponding vertices 
of the crystal graph (Proposition~\ref{SteinbergVertex}). Another consequence is 
stated in Corollary~\ref{ComparingDifferentCrystalGraphs}. 
Suppose that~$\lambda$ labels a weakly cuspidal module of~$G$ and that the 
$2$-core of~$\lambda$ is different from~$\lambda$ and contains more than one
node. Then there is a 
particular $e$-hook of~$\lambda$ such that the partition~$\lambda'$ obtained 
from~$\lambda$ by removing this $e$-hook also labels a weakly cuspidal module, 
and the two weakly cuspidal modules should give rise to isomorphic 
Harish-Chandra branching graphs. This is remarkable as~$n$ and~$n-e$ have 
different parities and the modules of $G = \GU_n(q)$ and $\GU_{n-e}(q)$ are not 
directly related via Harish-Chandra induction. We prove in 
Theorem~\ref{ComparingDifferentCrystalGraphs1} that, as predicted 
in~\ref{DifferentGraphs}, the two connected components in question are 
isomorphic (as unlabelled) graphs. A further consequence of our conjectures
is stated in Corollary~\ref{DifferentBlocks}: non-isomorphic composition 
factors of the socles of modules Harish-Chandra induced from $G = \GU_n(q)$ to 
$\GU_{n+2}(q)$, lie in different $\ell$-blocks. 

Let us finally comment on the history of this paper. First notes of the second
author date back to 1993, following the completion of~\cite{GeHiMaI}. There,
a general conjecture for the distribution of the simple modules of a unitary 
group into Harish-Chandra series for the linear prime case was presented. This 
conjecture was later verified in \cite{GruHi}. A further conjecture 
of~\cite{GeHiMaI} for the case that~$\ell$ divides $q + 1$ was proved 
in~\cite{GeHiMaII}. The conjectures in~\cite{GeHiMaI} were based on explicit 
decomposition matrices of unipotent modules of $\GU_n(q)$, computed by Gunter 
Malle. These decomposition matrices were completely known in the linear prime 
case for $n \leq 10$ and published in~\cite{GeHiMaI}. At that time, the 
information in the non-linear prime case was less comprehensive. Much more 
complete versions of these decomposition matrices and the distribution of
the unipotent modules into Harish-Chandra series are now available by the
recent work~\cite{DuMa} of Dudas and Malle.

Since the publication 
of~\cite{GeHiMaI}, many attempts have been made to find the combinatorial 
pattern behind the Harish-Chandra series of the unitary groups. The breakthrough
occurred in 2009, when the second and last author shared an office during a special 
program at the Isaac Newton Institute in Cambridge. The paper \cite{GeckPS} by
Geck and some other considerations of the second author suggested that the
simple modules of certain Iwahori-Hecke algebras of type~$B$ should label
some unipotent modules of the unitary groups. The paper~\cite{GeJa} by Geck 
and the third author on canonical basic sets then proposed the correct 
labelling by Uglov bipartitions. This set of bipartitions is defined 
through a certain crystal graph, called $\mathcal{G}_{\mathbf{c},e}$
below. The two authors compared their results on these crystal graphs
on the one hand, and on the known Harish-Chandra distribution on the other 
hand. Amazingly, the two results matched.

\section{A generalization of Harish-Chandra theory}
\label{HCGeneralization}

Here we introduce a generalization of Harish-Chandra theory for certain
families of classical groups by restricting the set of Levi subgroups. 

\subsection{}
\label{TheGroups}

Let~$q$ be a power of the prime~$p$. For a non-negative integer~$n$ let 
$G := G_n := G_n(q)$ denote one of the following classical groups, where we 
label the cases according to the (twisted) Dynkin type of the groups: 
\begin{description}
\item[$({^2\!A}_{2n-1})$] $\GU_{2n}(q)$, 
\item[$({^2\!A}_{2n})$] $\GU_{2n+1}(q)$, 
\item[$(B_{n})$] $\SO_{2n+1}(q)$, 
\item[$(C_{n})$] $\Sp_{2n}(q)$.
\end{description}
(We interpret $\GU_{0}(q)$ and $\Sp_{0}(q)$ as the trivial group.)

If $n \geq 1$, the group~$G$ is a finite group with a split $BN$-pair of 
characteristic~$p$, satisfying the commutator relations. In these cases, the 
Weyl group~$W$ of~$G$ is a Coxeter group of type~$B_n$, and we number the set 
$S = \{s_1, \ldots , s_n \}$ of fundamental reflections of~$W$ according to 
the following scheme.
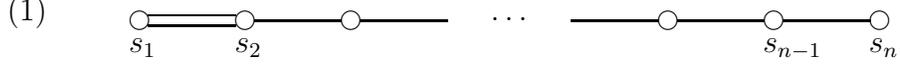
\begin{equation}
\label{TypeB}
\begin{picture}(260,20)(0,0)
\put(   0,   0){\circle{7}}                                                     
\put(  40,   0){\circle{7}}
\put(  80,   0){\circle{7}}
\put( 200,   0){\circle{7}}
\put( 240,   0){\circle{7}}                                   
\put( 280,   0){\circle{7}}                                   
\put( 3.5,   2){\line(1,0){ 33}}                                     
\put( 3.5,  -2){\line(1,0){ 33}}                         
\put(43.5,   0){\line(1,0){ 33}}                             
\put(83.5,   0){\line(1,0){ 33}}                             
\put( 129,   0){\mbox{\ $\ldots$\ }}
\put(163.5,  0){\line(1,0){ 33}}
\put(203.5,  0){\line(1,0){ 33}}
\put(243.5,  0){\line(1,0){ 33}}
\put(   -4,-12){$s_1$}
\put(   36,-12){$s_2$}
\put(  236,-12){$s_{n-1}$}                
\put(  276,-12){$s_{n}$}                
\end{picture}                                                         
\end{equation}

\bigskip
\bigskip

\subsection{}
\label{PureLevis}
A subset $I \subseteq S$ is called {\em left connected}, if it is of the form
$I = \{ s_1, s_2, \ldots , s_r \}$ for some $0 \leq r \leq n$. The corresponding
standard Levi subgroup~$L_I$ of~$G$ is denoted by $L_{r,n-r}$. A Levi 
subgroup~$L$ of~$G$ is called {\em pure}, if it is conjugate in~$N$ to a 
standard Levi subgroup~$L_I$ with~$I$ left connected. The set of all pure Levi 
subgroups of~$G$ is denoted by $\mathcal{L}^*$, whereas~$\mathcal{L}$ denotes 
the set of all $N$-conjugates of all standard Levi subgroups of~$G$. If $L \in 
\mathcal{L}^*$, a {\em pure Levi subgroup of}~$L$ is an element 
$M \in \mathcal{L}^*$ with $M \leq L$.

Notice that the set of $N$-conjugacy classes in $\mathcal{L}^*$ is linearly
ordered in the following sense. Let $L, M \in \mathcal{L}^*$. Then $|L| < |M|$
if and only if there is $x \in N$ such that ${^x\!L} \leq M$. In particular,
$|L| = |M|$ if and only if~$L$ and~$M$ are conjugate in~$N$.

Put~$\delta := 2$, if $G_n(q) = \GU_n(q)$, and $\delta := 1$, otherwise. Then
the standard Levi subgroup $L_{r,n-r}$ of~$G$ has structure
$$L_{r,n-r} \cong G_r( q ) \times \GL_1( q^\delta ) \times \cdots \times \GL_1( q^\delta )$$
with $n - r$ factors $\GL_1( q^\delta )$, and with a natural embedding of the
direct factors of~$L_{r,n-r}$ into~$G$.

\begin{lem*}
Let~$I$ and~$J$ be two left connected subsets of~$S$, and let~$x \in D_{IJ}$,
where $D_{IJ} \subseteq W$ denotes the set of distinguished double coset 
representatives with respect to the parabolic subgroups~$W_I$ and~$W_J$ of~$W$.
Then ${^x\!I} \cap J$ is left connected.
\end{lem*}
\begin{prf}
We identify~$W$ with the set of permutations~$\pi$ of $\{ \pm i \mid 
1 \leq i \leq n \}$ satisfying $\pi( -i ) = - \pi( i )$ for all 
$1 \leq i \leq n$. If $J = \emptyset$, there is nothing to prove. Thus assume 
that $J = \{ s_1, \ldots, s_r \}$ for some $1 \leq r \leq n$. Then~$W_J$ is the
stabilizer of the subset $\{ \pm i \mid 1 \leq i \leq r \}$ and all the 
singletons not in this set. It follows that ${^xW}_I \cap W_J$ is the
stabilizer of a set $\{ \pm i \mid i \in Z \}$ and all the singletons not in 
this set, where $Z \subseteq \{ 1, \ldots , r \}$.

On the other hand, if $J' := {^x\!I} \cap J$, then ${^xW}_I \cap W_J = W_{J'}$, as 
$x \in D_{IJ}$. This implies that~${J'}$ is left connected, as otherwise~$W_{J'}$
would not be a stabilizer as above.
\end{prf}

G\"otz Pfeiffer has informed us of a different proof of the above result, 
using the descent algebra of~$W$. Pfeiffer's proof also applies to Weyl groups
of type~$A$ and~$D$.
\begin{prp*}
Let~$L$,~$M$ be pure Levi subgroups of~$G$, and let $x \in N$.
Then ${^x\!L} \cap M$ is a pure Levi subgroup of~$G$.
\end{prp*}
\begin{prf}
We may assume that $L = L_I$ and $M = L_J$ for $I, J \subseteq S$ left
connected. As ${^x\!L}_I \cap L_J$ is conjugate in~$N$ to ${^y\!L}_I \cap L_J$,
where $y \in D_{IJ}$, we may also assume that $x \in D_{IJ}$. Then
${^x\!I} \cap J$ is left connected by the lemma. This completes the proof.
\end{prf}

\subsection{}
\label{WeakHCSeries}
Let~$k$ be a field of characteristic $\ell \neq p \geq 0$, 
such that~$k$ is a splitting field for all subgroups of~$G$. We write
$kG$-mod for the category of finite-dimensional $kG$-modules.
It is known that Harish-Chandra philosophy for $kG$ carries over
to the situation where $\mathcal{L}$ is replaced by $\mathcal{L}^*$.
The first ideas in this direction go back to Grabmeier's thesis \cite{Grabm}, 
who replaced Green correspondence in symmetric groups by a generalized
Green correspondence with respect to Young subgroups. Further developments
are due to Dipper and Fleischmann \cite{DiFl}. A comprehensive treatment 
including several new aspects can be found in \cite[Chapter~1]{CaEn}. 
The crucial ingredient in this generalization is Proposition~\ref{PureLevis}.

Let $L \in \mathcal{L}$. We write $R_L^G$ and ${^*\!R}_L^G$ for Harish-Chandra 
induction from $kL$-mod to $kG$-mod and Harish-Chandra restriction from $kG$-mod
to $kL$-mod, respectively. For $X \in kL$-mod we put
$$H_k(L,X) := \End_{kG}( R_L^G( X ) )$$
for the endomorphism algebra of~$R_L^G(X)$.

Let $X \in kG$-mod. We say that~$X$ is {\em weakly cuspidal},
if ${^*\!R}_L^G(X) = 0$ for all $G \neq L \in \mathcal{L}^*$.
A pair $(L,X)$ with $L \in \mathcal{L}^*$ and~$X$ a weakly cuspidal
simple $kL$-module is called a {\em weakly cuspidal pair}.
Let $(L,X)$ be a weakly cuspidal pair. Then the {\em weak Harish-Chandra
series} defined by $(L,X)$ consists of the simple $kG$-modules which are
isomorphic to submodules of $R_L^G(X)$. If $Y \in kG$-mod lies
in the weak Harish-Chandra series defined by $(L,X)$, then~$L \in \mathcal{L}^*$
is minimal with ${^*\!R}_L^G(Y) \neq 0$, and~$X$ is a composition factor of
${^*\!R}_L^G(Y)$.

We collect a few important facts about weak Harish-Chandra series.

\begin{prp*}
Let $(L,X)$ be a weakly cuspidal pair. 

{\rm (a)} Write
$$R_L^G( X ) = Y_1 \oplus \cdots \oplus Y_r$$
with indecomposable modules $Y_i$, $1 \leq i \leq r$. Then
each $Y_i$ has a simple head $Z_i$, which is also isomorphic to the socle
of~$Y_i$. Moreover, $Y_i \cong Y_j$, if and only if $Z_i \cong Z_j$.
The Harish-Chandra series defined by $(L,X)$ consists of the $kG$-modules
isomorphic to the $Z_i$.

{\rm (b)} The weak Harish-Chandra series partition the set of
isomorphism types of the simple $kG$-modules.

{\rm (c)} The weak Harish-Chandra series defined by $(L,X)$
is contained in a usual Harish-Chandra series, and thus every usual
Harish-Chandra series is partitioned into
weak Harish-Chandra series.
\end{prp*}
\begin{prf}
It follows from \cite[Theorems $1.20(iv)$, $2.27$]{CaEn} that $H_k(L,X)$ 
is a symmetric $k$-algebra (notice that the cited results are also valid
in our situation where $\mathcal{L}$ is replaced by~$\mathcal{L}^*$). 
This implies the statements of~(a) (see, e.g.\ \cite[Theorem~$1.28$]{CaEn}).

The proof of~(b) is analogous to the proof in the usual Harish-Chandra
theory.

To prove~(c), let $M \in \mathcal{L}$, and let $Z \in kM$-mod be 
cuspidal (in the usual sense) such that~$X$ occurs in 
the socle of $R_M^L( Z )$. Then $R_L^G( X )$ is a submodule of 
$R_L^G( R_M^L( Z ) ) \cong R_M^G( Z )$, and thus every simple module
in the socle of $R_L^G( X )$ also occurs in the socle of $R_M^G( Z )$
and hence in the usual Harish-Chandra series defined by $(M,Z)$.
\end{prf}

\subsection{}
\label{OtherCompositionFactors}
Let $(L,X)$ be a weakly cuspidal pair.
The following proposition gives information about those composition factors 
of $R_L^G(X)$ that do not lie in the weak Harish-Chandra series defined 
by $(L,X)$. The corresponding result for usual Harish-Chandra series is 
implicitly contained in \cite[Lemma~$5.7$]{HCseries} (see the remarks 
in~\cite[(2.2)]{GeHiMaI}). Since this result is particularly relevant in 
the definition of the Harish-Chandra branching graph, and since it is not 
explicitly formulated in \cite[Lemma~$5.7$]{HCseries}, and wrongly stated
in \cite[Proposition~$2.11(b)$]{GeHi}, we give a proof here.

\begin{prp*}
Let $(L,X)$ be a weakly cuspidal pair, and let~$Y$ be a composition
factor of $R_L^G(X)$. 
Suppose that~$Y$ lies in the weak Harish-Chandra series defined by 
$(M,Z)$, a weakly cuspidal pair.

Then there is $x \in N$ such that ${^x\!L} \leq M$. If ${^x\!L} = M$,
then $Z \cong {^x\!X}$. In particular, if~$Y$ does not lie in the weak
Harish-Chandra series defined by $(L,X)$, then $|L| < |M|$.
\end{prp*}
\begin{prf}
Let $P(Z)$ denote the projective cover of~$Z$. We have
$$0 \neq [ P(Z), {^*\!R}_M^G( Y ) ] = [ R_M^G ( P(Z) ), Y ],$$
the inequality arising from the fact that~$Z$ is a composition factor of 
${^*\!R}_M^G( Y )$, the equation arising from adjointness.
As~$Y$ is a composition factor of $R_L^G( X )$, we obtain
$$0 \neq [ R_M^G ( P(Z) ), R_L^G( X ) ] =
\sum_{x \in D_{M,L}} 
[ P( Z ), R^M_{M \cap {^x\!L}} ({^*\!R}_{M \cap {^x\!L}}^{^x\!L}( {^x\!X} ) ) ].$$
(Here, $D_{M,L} \subseteq N$ denotes a suitable set of representatives for 
double cosets with respect to parabolic subgroups of~$G$ with Levi 
complements~$M$ and~$L$, respectively.)
Thus there is $x \in D_{M,L}$ such that
$[ P( Z ), R^M_{M \cap {^x\!L}} ({^*\!R}_{M \cap {^x\!L}}^{^x\!L}( {^x\!X} ) ) ]
\neq 0$. As $(L,X)$ is a weakly cuspidal pair, so is $( {^x\!L}, {^x\!X} )$.
It follows that $M \cap {^x\!L} = {^x\!L}$, and thus ${^x\!L} \leq M$.
If ${^x\!L} = M$, we obtain $[ P(Z), {^x\!X} ] \neq 0$, hence our claim.
\end{prf}

\subsection{}
\label{AnExample}
If $\mbox{\rm char}(k) = 0$, a $kG$-module is {\em unipotent}, if it is simple 
and its character is unipotent. If $\ell > 0$, a $kG$-module is {\em unipotent}, 
if it is simple and its Brauer character (with respect to a suitable $\ell$-modular
system) is a linear combination of unipotent characters (restricted to $\ell'$-elements).

As $\mathcal{L}^* \subseteq \mathcal{L}$, every cuspidal $kG$-module~$X$ is
weakly cuspidal. The converse is not true, as the following example shows.
Let $G = \GU_6(q)$ and suppose that $\ell > 6$ and divides $q^2 - q + 1$.
The Levi subgroup $L = \GL_3( q^2 )$ (a Levi
complement of the stabilizer of a maximal isotropic subspace of the natural
vector space of~$G$), contains a cuspidal unipotent $kL$-module~$X$ by
\cite[Theorem 7.6]{GeHiMaI}. By applying \cite[Lemma~3.16]{GeHiMaII} and
\cite[Proposition 2.3.5]{GruGre} we find that $R_L^G( X )$ is indecomposable.
Let~$Y$ denote the unique head composition factor of~$R_L^G( X )$ (see
\cite[Theorem~$2.4$]{GeHiMaII}). By construction,~$Y$ is not cuspidal, but
weakly cuspidal. (The $kG$-module~$Y$ has label $2^3$ in the notation of
\cite[Table~$8$]{DuMa}).

Now suppose that $\ell = 0$. Then a weakly cuspidal unipotent $kG$-module is
cuspidal. Indeed, $\GL_n(q^\delta)$ has a cuspidal unipotent module over~$k$ 
only if $n = 1$. In particular, if $L \in \mathcal{L}$ has a cuspidal unipotent
module over~$k$, then $L \in \mathcal{L}^*$. If~$X$ is a weakly cuspidal 
unipotent $kG$-module and $L \in \mathcal{L}$ is minimal with ${^*\!R}_L^G(X) 
\neq 0$, every constituent of ${^*\!R}_L^G(X)$ is cuspidal.  Thus $L \in 
\mathcal{L}^*$ and hence, as $X$ is weakly cuspidal, $L = G$.

\section{The endomorphism algebra of Harish-Chandra induced 
weakly cuspidal modules}
\label{EndomorphismAlgebra}

In important special cases the endomorphism algebras $H_k(L,X)$ of weakly 
cuspidal pairs $(L,X)$ are Iwahori-Hecke algebras. The result applies in 
particular when~$X$ is unipotent. 

We keep the notation of Section~\ref{HCGeneralization},
except that we assume that $n \geq 1$ here. Thus if $G = G_n(q)$ is one of the 
groups introduced in~\ref{TheGroups}, then~$G$ has a split $BN$-pair of 
rank~$n$. Let~$\ell$ be a prime not dividing~$q$. We choose an $\ell$-modular
system $(K, \mathcal{O},k)$ such that~$K$ is large enough for~$G$. That is,
$\mathcal{O}$ is a complete discrete valuation ring with field of fractions~$K$
of characteristic~$0$, and residue class field~$k$ of characteristic~$\ell$.
Moreover,~$K$ is a splitting field for all subgroups of~$G$.

\subsection{}
\label{PreparingForQ}
Put $r := n - 1$ and $L := L_{r,1} \in \mathcal{L}^*$.
Thus $L = M \times T$ with $M \cong G_r(q)$ and $T \cong \GL_1(q^\delta)$.
(In case $n = 1$, either~$M$ is the trivial group, or cyclic of order $q + 1$ 
if $G = \GU_3(q)$.)
Let $P$ denote the standard parabolic subgroup of~$G$ with Levi complement~$L$ 
and let~$U$ denote its unipotent radical. We have $|W_G(L)| = 2$ and we let 
$s \in N_G(L)$ denote an inverse image of the involution in~$W_G(L)$. We 
choose~$s$ of order~$2$ if~$G$ is unitary or orthogonal, and of order~$4$ with
$s^2 \in T$ if~$G$ is symplectic, and such that~$s$ centralizes~$M$. (Such 
an~$s$ always exists.) 

Let~$R$ be one of the rings~$K$,~$\mathcal{O}$, or~$k$. As~$M$ is an epimorphic 
image of~$P$, we get a surjective homomorphism $\pi : RP \rightarrow RM$. 
Consider the element
\begin{equation}
\label{ComplicatedSum}
y := \sum_{\begin{array}{c} u,u' \in U\\su'sus \in P\end{array}} su'sus \in RP.
\end{equation}
Then $z := \pi(y) \in Z(RM)$ as~$s$ centralizes~$M$. 
\begin{lem*}
With the above notation, $z = (q-1)z'$ for some $z' \in Z(RM)$. In case~$G$ is 
a unitary group, we have $z' = 1 + (q + 1)z''$ for some $z'' \in Z(RM)$. 
\end{lem*}
\begin{prf}
We first claim that 
$T \cong \GL_1( q^\delta )$ acts on 
$$\mathcal{U} := \{ (u',u) \in U \times U \mid su'sus \in P \}$$ 
by 
$$x.(u',u) := (sxs^{-1}u'sx^{-1}s^{-1},xux^{-1}),
\quad x \in T, (u',u) \in \mathcal{U}.$$
Indeed,
\begin{equation}
\label{TActionOnU}
s(sxs^{-1}u'sx^{-1}s^{-1})s(xux^{-1})s = (s^2xs^{-2})su'sus(s^{-1}x^{-1}s)
\end{equation}
for $x \in T, (u',u) \in \mathcal{U}$.
As~$s$ normalizes~$T$, the claim 
follows. Now $\pi( x ) = 1$ for $x \in T$ and thus~(\ref{TActionOnU}) implies
$\pi( su'sus ) = \pi( sv'sv )$ if $(u',u), (v',v) \in \mathcal{U}$ are in the 
same $T$-orbit. 

The claims in the arguments below can be verified by a direct computation 
in~$G$. Suppose that~$G$ is a unitary or symplectic group. For each
$1 \neq u \in Z(U)$ there is a unique $u' \in Z(U)$ such that $(u',u) \in
\mathcal{U}$. For every such pair we have $\pi(su'sus) = 1$. The elements 
$(u',u) \in \mathcal{U}$ with $u \not\in Z(U)$ lie in regular $T$-orbits, as~$T$
acts fixed point freely on $U \setminus Z(U)$ by conjugation. This 
implies our result, as $|Z(U)| = q$ and $|T| = q^\delta - 1$. Now suppose 
that~$G$ is an orthogonal group. Then~$T$ acts with regular orbits 
on $U \setminus \{ 1 \}$, hence on $\mathcal{U}$, again implying our result. 
\end{prf}

\subsection{}
\label{StructureHLX}
Let $R$ be one of $K$ or~$k$. If~$X$ is an indecomposable $RG$-module, we
let~$\omega_X$ denote the central character of~$RG$ determined by the block 
containing~$X$.

Let~$r$ be an integer with $0 \leq r \leq n$ and put 
$m := n-r$. Let $L := L_{r,m} \in \mathcal{L}^*$ denote the standard Levi 
subgroup of $G = G_n(q)$ isomorphic to $G_r(q) \times \GL_1(q^\delta)^m$. 
Write~$M$ and~$T$ for the direct factors of~$L$ isomorphic to~$G_r(q)$ and
$\GL_1(q^\delta)^m$, respectively. Let~$X$ be a 
weakly cuspidal simple $RM$-module, extended trivially to an $RL$-module. 

For $R = K$ and~$X$ cuspidal, the following result is due to Lusztig (see 
\cite[Section~$5$]{LuszClass}).
\begin{thm*}
With the above notation,
$H_R(L,X)$ is an Iwahori-Hecke algebra corresponding to the Coxeter group 
of type~$B_m$, with parameters as in the following diagram.
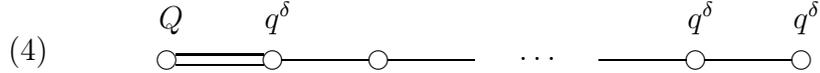
\begin{equation}
\label{Parameters}
\begin{picture}(240,30)(0,0)
\put(   0,   0){\circle{7}}                                                     
\put(  40,   0){\circle{7}}
\put(  80,   0){\circle{7}}
\put( 200,   0){\circle{7}}
\put( 240,   0){\circle{7}}                                   
\put( 3.5,   2){\line(1,0){ 33}}                                     
\put( 3.5,  -2){\line(1,0){ 33}}                         
\put(43.5,   0){\line(1,0){ 33}}                             
\put(83.5,   0){\line(1,0){ 33}}                             
\put( 129,   0){\mbox{\ $\ldots$\ }}
\put(163.5,  0){\line(1,0){ 33}}
\put(203.5,  0){\line(1,0){ 33}}
\put(   -3,+12){$Q$}
\put(   37,+12){$q^\delta$}
\put(  197,+12){$q^\delta$}                                               
\put(  237,+12){$q^\delta$}                
\end{picture}                                                         
\end{equation}

\bigskip

\noindent
The parameter~$Q$ is determined as follows. Let~$U$ and~$z$ be as 
in~\ref{PreparingForQ}, applied to $G_{r+1}$. Put $\gamma := \omega_X( z ) 
\in R$ and let $\xi \in R$ be a solution of the quadratic equation
$$x^2 - \gamma\,x - |U| = 0.$$
Then 
$$Q = \frac{\xi\gamma}{|U|} + 1.$$
Moreover, the following statements hold.

{\rm (a)} Suppose that $R = k$ and that $X$ lies in a block containing a
cuspidal $KM$-module~$Y$. If $\hat{Q}$ is the parameter of
$H_K(L,Y)$ associated to the leftmost node of the diagram~(\ref{Parameters}),
then $Q$ is the reduction modulo~$\ell$ of~$\hat{Q}$.

{\rm (b)} If $R = k$ and $\ell \mid q - 1$, then $Q = 1$.

{\rm (c)} If $R = k$ and $\ell \mid q + 1$, then $Q = -1$.
\end{thm*}
\begin{prf}
First notice that we have $W_G( L, X ) = W_G( L )$, and that $W_G(L)$ is 
isomorphic to a subgroup of~$W$ and a Coxeter group of type~$B_m$ (see 
\cite{How}). We also have 
$$\mbox{\rm dim}_R( H_R(L,X) ) = |W_G( L )|.$$
Put $\mathcal{N}(L) := (N_G(L) \cap N)L$ (recall that $G$ has a $BN$-pair),
so that $W_G(L) = \mathcal{N}(L)/L$. Then $\mathcal{N}(L) = M \times C$
with $T \leq C$ and $C/T \cong W_G(L)$. In particular, we may view~$X$ as an
$R\mathcal{N}(L)$-module on which~$C$ acts trivially.

The parameters not corresponding to the leftmost node of~(\ref{Parameters})
can now be computed exactly as in the case where~$X$ is cuspidal and unipotent 
(see \cite[Proposition~$4.4$]{GeHiMaII}). 

To determine~$Q$ we may assume that $m = 1$. Thus $G = G_{r+1}(q)$ and
$L \cong M \times \GL_1(q^2)$. We are thus in the situation 
of~\ref{PreparingForQ} and make use of the notation introduced there. 
Then $H := H_R(L,X)$ is $2$-dimensional over~$R$ with 
basis elements~$B_1$ and~$B_s$, where~$B_1$ is the unit element of~$H$ 
and~$B_s$ is defined as follows. We may realize~$R_L^G(X)$ as
$$R_L^G(X) = \{ f : G \rightarrow X \mid f(hg) = h.f(g),
\mbox{\rm\ for all\ } h \in P, g \in G \}.$$
Then $B_{s}$ is defined by
$$B_{s}(f)(g) := \frac{1}{|U|}\sum_{u \in U}f(sug), \quad\quad 
f \in R_L^G(X), g \in G,$$
as $s \in C$ acts trivially on~$X$.
We have $B_{s}^2 = 
\zeta B_1 + \eta B_{s}$ with $\zeta = 1/|U|$, and~$\eta$ such that
the element~$y$ of~(\ref{ComplicatedSum}) 
acts as the scalar $|U| \eta$ on~$X$. This is proved exactly as in
\cite[Proposition 3.14]{HowLeh}. 

Now~$y$ acts in the same way on~$X$ as $z = \pi(y)$.
Since~$X$ is absolutely irreducible, $z \in Z(RM)$ acts by the scalar
$\omega_X( z )$. Thus $|U|\eta = \omega_X( z ) = \gamma$.
Put
$$T_s := \xi B_s, \quad T_1 := B_1.$$
Then
$$T_s^2 = Q T_1 + (Q - 1)T_s$$
with $Q = \xi\eta + 1$. This gives our first claim.

To prove~(a), put $\hat{\gamma} := \omega_{Y}( z )$, and let
$\hat{\xi}$ be a solution of 
$x^2 - \hat{\gamma} x - |U| = 0$.
Observe that $\hat{\gamma}, \hat{\xi} \in \mathcal{O}$. Then the
reduction modulo $\ell$ of $\hat{\gamma}$ equals $\omega_X( z )$,
and the reduction modulo $\ell$ of $\hat{\xi}$ is a solution of
$x^2 - \gamma x - |U| = 0$. Thus the reduction modulo $\ell$ of
$\hat{Q} := \hat{\xi}\hat{\eta} + 1$ equals $\xi\eta + 1 = Q$
and~(a) is proved.

Suppose now that $R = k$. If $\ell \mid q - 1$, we have $\gamma = 0$ by 
Lemma~\ref{PreparingForQ} and thus $Q = 1$. If~$G$ is unitary and $\ell 
\mid q + 1$, we have $\gamma = -2$, again by Lemma~\ref{PreparingForQ}. 
Also, $|U|$ is an odd power of~$q$, i.e.\ $|U| = - 1$ in~$k$, hence 
$\xi = - 1$ and $Q = - 1$. This completes our proof.
\end{prf}

\section{The Harish-Chandra branching graph}

In this section we fix a prime power~$q$ of~$p$ and a prime $\ell \neq p$.
We also let~$k$ denote an algebraically closed field of characteristic~$\ell$.

\subsection{}
For $n \in \mathbb{N}$, we let $G := G_n := G_n(q)$ denote one of the 
groups of~\ref{TheGroups}. Recall that $G_n$
is naturally embedded into $G_{n+1}$, by embedding $G_n$ into the
pure Levi subgroup $L_{n,1} \cong G_n \times \GL_1( q^\delta )$ 
of~$G_{n+1}$. By iterating, we obtain an embedding of $G_n$ into
$G_{n+m}$ for every $m \in \mathbb{N}$.

By $kG\mbox{\rm -mod}^u$ we denote the full subcategory of 
$kG\mbox{\rm -mod}$ consisting of the modules that have a filtration by
unipotent $kG$-modules. By the result of Brou{\'e} and Michel \cite{BM}, 
and by \cite{Hi}, $kG\mbox{\rm -mod}^u$ is a direct sum of blocks of~$kG$.
The above embedding of $G_n$ into $G_{n+m}$ yields a functor
$$R_n^{n+m}: kG_n\mbox{\rm -mod}^u \rightarrow kG_{n+m}\mbox{\rm -mod}^u,$$
defined by 
$$R_n^{n+m}( X ) := R_{L_{n,m}}^{G_{n+m}}( \Infl_{G_n}^{L_{n,m}}( X )),\quad
X \in kG_n\mbox{\rm -mod}^u,$$
where $\Infl_{G_n}^{L_{n,m}}( X )$ denotes the trivial extension of~$X$ to
$L_{n,m} \cong G_n \times \GL_1( q^\delta )^m$.
The adjoint functor 
$${^*\!R}_n^{n+m}: kG_{n+m}\mbox{\rm -mod}^u \rightarrow kG_{n}\mbox{\rm -mod}^u,$$
is given by
$$R_n^{n+m}( X ) := \Res_{G_n}^{L_{n,m}}( {^*\!R}_{L_{n,m}}^{G_{n+m}}( X ) ),\quad
X \in kG_{n+m}\mbox{\rm -mod}^u.$$

Let $\mathcal{R}_n := \mathcal{R}_n(q)$ denote the Grothendieck group of 
$kG_n\mbox{\rm -mod}^u$, and put 
$$\mathcal{R} := \mathcal{R}(q) := \bigoplus_{n \in \mathbb{N}} \mathcal{R}_n.$$
For an object $X \in kG_n\mbox{\rm -mod}^u$, we let $[X]$ denote its
image in $\mathcal{R}_n$. 

\subsection{}
\label{HCBG}
The {\em (twisted) Dynkin type} of~$G$ is one of the symbols ${^2\!A}_\iota$
with $\iota \in \{ 0, 1 \}$, $B$ or $C$, where $\GU_r(q)$
has twisted Dynkin type ${^2\!A}_\iota$ with $\iota \equiv (r \mbox{\rm\ mod\ } 2)$.

The {\em Harish-Chandra branching graph} $\mathcal{G}_{\mathcal{D},q,\ell}$ 
corresponding to~$q$,~$\ell$ and the (twisted) Dynkin type~$\mathcal{D}$
is the directed graph whose vertices are the elements~$[X]$, where~$X$ is a 
simple object in $kG_n\mbox{\rm -mod}^u$ for some $n \in \mathbb{N}$. Thus the 
vertices of $\mathcal{G}_{\mathcal{D},q,\ell}$ are the standard basis elements 
of~$\mathcal{R}$. We say that the a vertex $[X]$ has {\em rank}~$n$, if $[X] 
\in \mathcal{R}_n$. 
Let $[X]$ and~$[Y]$ be vertices in $\mathcal{G}_{\mathcal{D},q,\ell}$. Then
there is a directed edge from $[X]$ to~$[Y]$ if and only if there is $n \in 
\mathbb{N}$ such that~$[X]$ has rank~$n$ and~$[Y]$ has rank $n+1$, and such 
that~$Y$ is a head composition factor of $R_n^{n+1}( X )$. 
A vertex in $\mathcal{G}_{\mathcal{D},q,\ell}$ is called a {\em source vertex}, 
if it has only outgoing edges.

As every unipotent $kG$-module is self dual,~$Y$ 
is a head composition factor of $R_n^{n+1}( X )$ if and only if~$Y$ is in the 
socle of $R_n^{n+1}( X )$. By adjunction,~$Y$ is a head composition factor of 
$R_n^{n+1}( X )$ if and only if~$X$ is in the socle of ${^*\!R}_n^{n+1}( Y )$, 
and~$Y$ is in the socle of $R_n^{n+1}( X )$ if and only if~$X$ is a head 
composition factor of ${^*\!R}_n^{n+1}( Y )$. 

An example for part of a Harish-Chandra branching graph is displayed in
Table~\ref{BranchingGraphExample}, where the vertices are represented by 
their labels. This can be proved with the help of the decomposition matrices 
computed in~\cite{DuMa} plus some ad hoc arguments.

\subsection{}
\label{PathsInHCGraph}
We have the following relation with the weak Harish-Chandra series of~$G$.
\begin{prp*}
Let~$[X]$ be a vertex of rank~$n$ of $\mathcal{G}_{\mathcal{D},q,\ell}$.
Then~$[X]$ is a source vertex if and only if $X \in kG_n\mbox{\rm -mod}^u$
is weakly cuspidal. 

Suppose that~$X$ is weakly cuspidal and let $m \in \mathbb{N}$. View~$X$ as a 
module of $L_{n,m}$ via inflation. Then a simple object 
$Y \in kG_{n+m}\mbox{\rm -mod}^u$
lies in the weak $(L_{n,m}, X)$ Harish-Chandra series, if and only if there
is a directed path from~$[X]$ to~$[Y]$ in $\mathcal{G}_{\mathcal{D},q,\ell}$. 
\end{prp*}
\begin{prf}
Clearly,~$X$ is weakly cuspidal if $n = 0$. Assume that $n \geq 1$. Then~$X$ is 
weakly cuspidal if and only if ${^*\!R}_{n-1}^{n}( X ) = 0$, which is the case
if and only if~$[X]$ is a source vertex.

Assume now that~$X$ is weakly cuspidal, let $m \in \mathbb{N}$ and let~$[Y]$
be a vertex of rank $n + m$. Suppose there is a path from~$[X]$ to~$[Y]$. 
We proceed by induction on~$m$ to show that~$Y$ occurs in the head of $R_{n}^{n+m}( X )$. 
If $m = 0$, there is nothing to prove. So assume that $m > 0$ and that the 
claim has been prove for~$m - 1$. Let $[Z]$ be a vertex of rank $n + m - 1$ that
occurs in a path from~$[X]$ to~$[Y]$. By induction,~$Z$ is a head composition
factor of $R_{n}^{n+m-1}( X )$. By exactness, $R_{n+m-1}^{n+m}( Z )$ is
a quotient of $R_{n+m-1}^{n+m}( R_{n}^{n+m-1}( X ) ) \cong R_{n}^{n+m}( X )$.
As~$Y$ is a quotient of $R_{n+m-1}^{n+m}( Z )$, we are done.

Suppose now that~$Y$ occurs in the head of $R_{n}^{n+m}( X )$. We proceed 
by induction on~$m$ to show that there is a path from~$[X]$ to~$[Y]$, the
cases $m \leq 1$ being trivial. As~$Y$ is isomorphic to a quotient of
$R_{n}^{n+m}( X ) \cong R_{n+m-1}^{n+m}( R_{n}^{n+m-1}( X ) )$, there is
a composition factor~$Z$ of $R_{n}^{n+m-1}( X )$ such that~$Y$ is a quotient
of~$R_{n+m-1}^{n+m}( Z )$. In particular, there is an edge from~$[Z]$ to~$[Y]$.
If~$Z$ occurs in the head of~$R_{n}^{n+m-1}( X )$, there is a path from~$[X]$ 
to~$[Z]$ by induction, and we are done. Aiming at a contradiction, assume 
that~$Z$ does not occur in the head of~$R_{n}^{n+m-1}( X )$. Then~$Z$ does
not lie in the weak Harish-Chandra series of $G_{n+m-1}$ defined by
$(L_{n,m-1}, X )$. It follows from Proposition~\ref{OtherCompositionFactors}
that~$Z$ lies in the weak Harish-Chandra series defined by 
$(L_{n',n-n' + m - 1}, X')$ for some $n < n'$ and some weakly cuspidal 
module~$X'$. In particular,~$Y$ lies in this weak Harish-Chandra series. 
This contradiction completes our proof.
\end{prf}

\section{Conjectures}
\label{Conjectures}

Here we formulate a series of conjectures about the $\ell$-modular Harish-Chandra
series and the Harish-Chandra branching graph for the unitary groups.

\subsection{}
\label{UnitaryGroups}
As always, we let~$q$ denote a power of a prime~$p$, and we fix a prime~$\ell$
different from~$p$. The multiplicative order of~$-q$ modulo~$\ell$ is denoted
by $e := e(q,\ell)$. Thus~$e$ is the smallest positive integer such that~$\ell$
divides $(-q)^e - 1$.

For a non-negative integer~$n$ we let $G := \GU_n( q )$ be the unitary group of 
dimension~$n$. Also, $(K, \mathcal{O}, k)$ denotes an $\ell$-modular system such
that~$K$ is large enough for~$G$ and with~$k$ algebraically closed.

\subsection{} \label{CombinatorialNotions} 
The set of partitions of a non-negative integer~$n$ is 
denoted by~$\mathcal{P}_n$ and we write $\lambda \vdash n$ if $\lambda \in 
\mathcal{P}_n$. We put $\mathcal{P} := \cup_{n \in \mathbb{N}} \mathcal{P}_n$.
Let $\lambda \in \mathcal{P}$. Then $\lambda_{(2)}$ and $\lambda^{(2)}$ denote
the $2$-core and the $2$-quotient of~$\lambda$, respectively. (As in 
\cite[Section~$1$]{fosri3}, the $2$-quotient is determined via a $\beta$-set 
for $\lambda$ with an odd number of elements, where we use the term $\beta$-set
in its original sense of being a finite set of non-negative integers as 
introduced in \cite[p.~$77$f]{jake}.)
For a non-negative integer~$t$ we write $\Delta_t := (t, t - 1, \ldots , 1)$ 
for the triangular partition of $t(t+1)/2$. Then $\lambda_{(2)} = \Delta_t$
for some $t \in \mathbb{N}$. Suppose that $\lambda^{(2)} = (\mu^1,\mu^2)$.
We then put $\bar{\lambda}^{(2)} := (\mu^1,\mu^2)$ if $t$ is even, and
$\bar{\lambda}^{(2)} := (\mu^2,\mu^1)$, otherwise.
If $\mu = (\mu^1,\mu^2)$ is a bipartition,
we let $\Phi_t( \mu )$ denote the unique partition $\lambda$ with 
$\lambda_{(2)} = \Delta_t$ and $\bar{\lambda}^{(2)} = (\mu^1,\mu^2)$
(see \cite[Theorem~$2.7.30$]{jake}).

The set of bipartitions of~$n$ is denoted by 
$\mathcal{P}_n^{(2)}$, and we put $\mathcal{P}^{(2)} := 
\cup_{n \in \mathbb{N}} \mathcal{P}_n^{(2)}$. Finally, we write 
$\mu \vdash_2 n$ if $\mu \in \mathcal{P}_n^{(2)}$.

\subsection{}
\label{FongSrinivasan}
By a result of Lusztig and Srinivasan \cite{LuSri}, the unipotent $KG$-modules 
are labelled by partitions of~$n$. We write $Y_\lambda$ for the unipotent 
$KG$-module labelled by $\lambda \in \mathcal{P}_n$. Let~$\lambda$ and $\mu$ 
be partitions of~$n$. It follows from the main result of Fong and Srinivasan
\cite[Theorem (7A)]{fosri1}, that $Y_\lambda$ and~$Y_\mu$ lie in the same 
$\ell$-block of~$G$, if and only if~$\lambda$ and~$\mu$ have the same $e$-core. 
The $e$-{\em weight} and the $e$-{\em core} of the $\ell$-block 
containing~$Y_\lambda$ are, by definition,  the $e$-weight and the $e$-core 
of~$\lambda$, respectively.

It was shown by Geck in \cite{TriangularShape} that if the~$Y_\lambda$, 
$\lambda \vdash n$, are ordered downwards lexicographically, the corresponding 
matrix of $\ell$-decomposition numbers is square and upper unitriangular. This 
defines a 
labelling of the unipotent $kG$-modules by partitions of~$n$, and we write 
$X_\mu$ for the unipotent $kG$-module labelled by $\mu \in \mathcal{P}_n$. 
Thus $X_\mu$ is determined by the following two conditions. Firstly,~$X_\mu$ 
occurs exactly once as a composition factor in a reduction modulo~$\ell$ 
of~$Y_\mu$, and secondly, if $X_\mu$ is a composition factor in a reduction 
modulo~$\ell$ of $Y_\nu$ for some $\nu \in \mathcal{P}_n$, then $\nu \leq \mu$.

\subsection{} \label{Compatibility} Our first conjecture asserts a compatibility
between ordinary and modular Harish-Chandra series.

\begin{conj*}
Let $\mu, \nu \in \mathcal{P}_n$.
If $X_\mu$ and $X_\nu$ lie in the same weak Harish-Chandra series of 
$kG$-modules, then $\mu$ and~$\nu$ have the same $2$-core, i.e.\ $Y_\mu$ and 
$Y_\nu$ lie in the same Harish-Chandra series of $KG$-modules. (In other words, 
the partition of $\mathcal{P}_n$ arising from the weak $\ell$-modular 
Harish-Chandra series is a refinement of the partition of $\mathcal{P}_n$ 
arising from the ordinary Harish-Chandra series.)
\end{conj*}

\subsection{}
\label{CuspidalityAndECore}
We also conjecture that a weakly cuspidal unipotent module can only occur in an 
$\ell$-block of~$G$ which contains a cuspidal simple $KG$-module (not 
necessarily unipotent). In fact, if~$e$ is odd, a unipotent 
$\ell$-block~$\mathbf{B}$ contains a cuspidal simple $KG$-module if and only 
if the $e$-core of~$\mathbf{B}$ is a $2$-core. This can be seen as follows.
Suppose first that the $e$-core of~$\mathbf{B}$ is the $2$-core
$\Delta_s$.
Put $m' := s(s+1)/2$. Let~$x$ be an $\ell$-element in~$G$ with $C := C_{G}( x ) 
= (q^e + 1)^w \times \GU_{m'}(q)$, where $(q^e + 1)^w$ denotes a direct product
of~$w$ factors of the cyclic group of order $q^e + 1$ (and $n = we + m'$).
Let~$Z$ denote the cuspidal unipotent $KC$-module labelled by~$\Delta_s$,
and let~$Y$ be the simple $KG$-module corresponding to~$Z$ under Lusztig's
Jordan decomposition. Then~$Y$ is cuspidal by~\cite[7.8.2]{LuszClass}, and~$Y$ 
lies in~$\mathbf{B}$ by~\cite[Theorem (7A) and Proposition (4F)]{fosri1}.
Conversely, suppose that~$\mathbf{B}$ contains some cuspidal simple 
$KG$-module~$Y$.
Then~$Y$ determines a unipotent $KC$-module, where~$C$ is the centralizer in~$G$
of some $\ell$-element. Let $\mu \in \mathcal{P}$ be the partition labelling~$Z$.
Then~$\mu$ is a $2$-core, and in turn, the $e$-core of~$\mu$ is a $2$-core
as well. As the $e$-core of~$\mu$ equals the $e$-core of~$\mathbf{B}$,
again by~\cite[Theorem (7A) and Proposition (4F)]{fosri1}, our claim follows.

\begin{conj*}
Let $\lambda \in \mathcal{P}_n$. If $X_\lambda$ is weakly cuspidal, then
the $e$-core of $\lambda$ is a $2$-core.
\end{conj*}

\noindent
It follows from~\cite[Corollary~$8.8$]{GruHi} that if $e$ is even, 
then~$X_\lambda$ is cuspidal if and only if $\lambda$ is a $2$-core.
(In this case,~$\lambda$ also is an $e$-core.)

Assuming that Conjecture~\ref{CuspidalityAndECore} holds, the parameter~$Q$
of a weakly cuspidal unipotent $kG$-module~$X_\lambda$ of~$G$ can be computed 
from the $e$-core of~$\lambda$ by Corollary~\ref{ParameterOf2Core} below.

\subsection{}
\label{Fockspace}
To present our next conjectures, we first have to introduce the Fock space of 
level~$2$ and its corresponding crystal graph. The results summarized below are 
due to Jimbo, Misra, Miwa and Okado \cite{JMMO} and Uglov \cite{Uglov}. For a 
detailed exposition see also \cite[Chapter~$6$]{GeJaHecke}.

A \textit{charged bipartition} is a pair $(\mu,\mathbf{c})$, written as 
$|\mu, \mathbf{c}\rangle$ with $\mu \in \mathcal{P}^{(2)}$ and $\mathbf{c} 
\in \mathbb{Z}^2$. Fix $\mathbf{c} = (c_1,c_2) \in \mathbb{Z}^2$, 
and let~$v$ denote an indeterminate. The Fock space (of level~$2$) and 
charge~$\mathbf{c}$ is the $\mathbb{Q}(v)$-vector space
$$\mathcal{F}_{\mathbf{c}} := \bigoplus_{m \in \mathbb{N}} 
\bigoplus_{\mu \vdash_2 m} \mathbb{Q}(v)|\mu,\mathbf{c}\rangle.$$

Assume that $e \geq 2$. There is an action of the quantum group 
$\mathcal{U}'_v( \widehat{\mathfrak{sl}_e} )$ on $\mathcal{F}_{\mathbf{c}}$
such that $\mathcal{F}_{\mathbf{c}}$ is an integrable 
$\mathcal{U}'_v( \widehat{\mathfrak{sl}_e} )$-module and
$|\mu,\mathbf{c}\rangle$ is a weight vector for every $m \in \mathbb{N}$
and $\mu \vdash_2 m$. Moreover, $|(-,-), \mathbf{c}\rangle$ is a highest 
weight vector and 
$\mathcal{U}'_v( \widehat{\mathfrak{sl}_e} ).|(-,-),\mathbf{c}\rangle$ 
is isomorphic to $V( \Lambda( \mathbf{c} ) )$, the simple highest weight module 
with weight $\Lambda( \mathbf{c} ) = \Lambda_{c_1 \text{\ mod\ e}} + 
\Lambda_{c_2 \text{\ mod\ e}}$, where the $\Lambda_i$ $0 \leq i \leq e - 1$ 
denote the fundamental weights of~$\widehat{\mathfrak{sl}_e}$. We write
$\mathcal{F}_{\mathbf{c},e}$ when we view $\mathcal{F}_{\mathbf{c}}$ as
a $\mathcal{U}'_v( \widehat{\mathfrak{sl}_e} )$-module.

There is a crystal graph $\mathcal{G}_{\mathbf{c},e}$ describing the 
canonical basis of $\mathcal{F}_{\mathbf{c},e}$. The vertices of 
$\mathcal{G}_{\mathbf{c},e}$ are all charged bipartitions
$|\mu,\mathbf{c}\rangle$, $\mu \vdash_2 m$, $m \in \mathbb{N}$. There is a 
directed, coloured edge 
$|\mu,\mathbf{c}\rangle \stackrel{i}{\rightarrow} |\nu,\mathbf{c} \rangle$
if and only if $\nu$ is obtained from~$\mu$ by adding a good $i$-node,
where the colours~$i$ are in the range $0 \leq i \leq e - 1$. The
associated Kashiwara operator $\tilde{f}_i$ acts on $\mathcal{G}_{\bc,e}$
by mapping the vertex $|\mu, \mathbf{c}\rangle$ to $|\nu, \mathbf{c}\rangle$
if and only if there is an edge 
$|\mu,\mathbf{c}\rangle \stackrel{i}{\rightarrow} |\nu,\mathbf{c} \rangle$,
and to~$0$, otherwise (see e.g.\ \cite[6.1]{GeJaHecke}).

Let us now describe, following~\cite{GeJaHecke}, how to compute the good 
$i$-nodes of $|\mu,\mathbf{c}\rangle$, and thus the 
graph~$\mathcal{G}_{\mathbf{c},e}$, algorithmically.
A \textit{node} of $\mu = (\mu^1,\mu^2)$ is a triple $(a,b,j)$, where $(a,b)$ 
is a node in the Young diagram of~$\mu^j$, for $j = 1, 2$. 
A node $\gamma$ of $\mu$ is called addable (respectively removable) if 
$\mu\cup\{\gamma\}$ (respectively $\mu\backslash\{\gamma\}$) is still a 
bipartition. The content of $\gamma = (a,b,j)$ is the integer $\cont(\gamma) =
b-a+c_j$. The residue of $\gamma$ is the element of $\{ 0, 1, \ldots , e - 1 \}$
defined by $\res(\gamma) = \cont(\gamma) 
\text{\rm\ mod\ } e$. For $0 \leq i \leq e-1$, $\gamma$ is called an $i$-node 
if $\res(\gamma)=i$.

Fix $i \in \{ 0, 1, \ldots , e - 1 \}$, and define an order on the set of 
addable and removable $i$-nodes of~$\mu$ by setting
$$\gamma \prec_\bc \gamma' \text{ if }   
\left\{      \begin{array}{l}      \cont(\gamma) < \cont(\gamma') \text{ or }
\\ \cont(\gamma) = \cont(\gamma') \mand j > j'.   \end{array}   \right.$$
Sort these set of nodes according to~$\prec_\bc$, starting from the smallest
one. Encode each addable (respectively removable) $i$-node 
by the letter $A$ (respectively $R$), and delete recursively all occurences of 
consecutive letters $RA$. This yields a word of the form $A^{\alpha_i}R^{\beta_i}$, 
which is called the reduced $i$-word of $\mu$. Note that by Kashiwara's crystal 
theory \cite[Section 4.2]{Kashiwara}, we have the following expression for 
the weight of the vector $|\mu,\bc\rangle$:
\begin{equation}\label{weight}\wt(\mu,\bc) = \sum_{i=0}^{e-1} (\alpha_i - \beta_i) 
\Lambda_i.\end{equation}
Let $\gamma$ be the rightmost addable (respectively leftmost removable) $i$-node 
in the reduced $i$-word of $\mu$. Then~$\gamma$ is called the \textit{good 
addable} (respectively \textit{good removable}) $i$-node of $\mu$. 

Each connected component of $\mathcal{G}_{\mathbf{c},e}$ is isomorphic to the 
crystal of a simple highest
weight module of $\mathcal{U}'_v( \widehat{\mathfrak{sl}_e} )$, whose highest
weight vector is the unique source vertex of the component.
The {\em rank} of a vertex $|\mu,\mathbf{c}\rangle$ of
$\mathcal{G}_{\mathbf{c},e}$ is~$m$, if $\mu \vdash_2 m$. We write
$\mathcal{G}_{\mathbf{c},e}^{\leq m}$ for the induced subgraph of
$\mathcal{G}_{\mathbf{c},e}$ containing the vertices of rank at most~$m$.

As an example, the graph $\mathcal{G}_{(0,0),3}^{\leq 3}$ is displayed in
Table~\ref{CrystalGraphExample}.

\subsection{}
\label{GraphComparison}
Let~$t$ be a non-negative integer, put $r := t(t+1)/2$ and $\iota := 
r\,(\mbox{\rm mod\ } 2) \in \{0,1\}$. Then $K\GU_r(q)$ has a unipotent 
cuspidal module~$Y$, and $( \GU_r(q), Y )$ determines a Harish-Chandra 
series of unipotent
$K\GU_{r + 2m}(q)$-modules for every $m \in \mathbb{N}$.
Recall from~\ref{HCBG} that $\mathcal{G}_{{^2\!A_\iota},q,\ell}$ denotes
the Harish-Chandra branching graph corresponding to~$q$, $\ell$ and the 
groups $\GU_{2n + \iota}(q)$.
As we are dealing exclusively with unitary groups in this section, we shall 
replace the index ${^2\!A}_\iota$ by~$\iota$ in the symbol for the graph. 
The vertices of
$\mathcal{G}_{\iota,q,\ell}$ correspond to the isomorphism classes of the
unipotent $k\GU_{2n + \iota}(q)$-modules, where~$n$ runs through the set of
positive integers. We may thus label the vertices of 
$\mathcal{G}_{\iota,q,\ell}$ by the set 
$\cup_{n \in \mathbb{N}} \mathcal{P}_{2n + \iota}$.

To formulate our next conjecture, we assume that Conjecture~\ref{Compatibility}
holds. Under this assumption, the induced subgraph of 
$\mathcal{G}_{\iota,q,\ell}$ whose vertices are labelled by the set of
partitions with $2$-core~$\Delta_t$, is a union of connected components
of $\mathcal{G}_{\iota,q,\ell}$. 
We write 
$\tilde{\mathcal{G}}_{\iota,q,\ell}^{t}$ for the graph with vertices 
$\mathcal{P}^{(2)}$, and a directed edge $\mu \rightarrow \nu$, if and only 
if there is a directed edge in $\mathcal{G}_{\iota,q,\ell}$ between the vertices 
labelled by $\Phi_t(\mu)$ and $\Phi_t(\nu)$. If $\mu \vdash_2 m$ is a vertex of
$\tilde{\mathcal{G}}_{\iota,q,\ell}^{t}$, the {\em rank} of this vertex is~$m$.
For a non-negative integer~$d$ we let 
$\tilde{\mathcal{G}}_{\iota,q,\ell}^{t,\leq d}$ denote the induced subgraph of
$\tilde{\mathcal{G}}_{\iota,q,\ell}^{t}$ containing the vertices of rank at 
most~$d$.

\begin{conj*}
Let the notation be as above. Assume that $e$ is odd and put $\mathbf{c} := 
(t + (1-e)/2,0)$. Then there is an integer $b := b(\ell)$ such that
$\tilde{\mathcal{G}}_{\iota,q,\ell}^{t,\leq b}$ equals 
$\mathcal{G}_{\mathbf{c},e}^{\leq b}$, if the colouring of the edges of
the latter graph is neglected.
\end{conj*}

\subsection{}
\label{Uglov}
As the Harish-Chandra series of unipotent $kG$-modules can be read off from
the Harish-Chandra branching graph by Proposition~\ref{PathsInHCGraph},
the truth of Conjecture~\ref{GraphComparison} would give an algorithm to 
determine the partition of the $kG$-modules into weak Harish-Chandra series
from the labels of the modules, at least if $\ell$ is large enough. In 
particular, the question of whether~$X_\lambda$ is weakly cuspidal, can be 
read off from~$\lambda$.

\begin{conj*}
Let $\lambda \in \mathcal{P}$ and let $t \in \mathbb{N}$ such that 
$\lambda_{(2)} = \Delta_t$. Let $\mu = \bar{\lambda}^{(2)}$ 
(see~\ref{CombinatorialNotions}).
Assume that $\ell$ is large enough, that~$e$ is odd and put $\mathbf{c} := 
(t + (1-e)/2,0)$.

Then $X_\lambda$ is weakly cuspidal, if and only if $|\mu,\mathbf{c}\rangle$
is a source vertex in $\mathcal{G}_{\mathbf{c},e}$.

Suppose that $X_\lambda$ is weakly cuspidal and let $\rho \in \mathcal{P}$. 
Then $X_\rho$ lies in the weak Harish-Chandra series defined by~$X_\lambda$,
if and only if $\rho_{(2)} = \lambda_{(2)} = \Delta_t$, and 
$|\bar{\rho}^{(2)},\mathbf{c}\rangle$ lies in the connected component of 
$\mathcal{G}_{\mathbf{c},e}$ containing $|\mu,\mathbf{c}\rangle$, i.e.\ 
$|\bar{\rho}^{(2)},\mathbf{c}\rangle$ is obtained from 
$|\mu,\mathbf{c}\rangle$ by adding a sequence of good nodes.
\end{conj*}

\section{Some evidence}
\label{Evidence}

Here we present the evidence for our conjectures. Keep the notation of 
Section~\ref{Conjectures}. We also assume that $e$ is odd and larger than~$1$
in this section.

\subsection{}
Conjecture~\ref{GraphComparison} holds for $e = 3, 5$ and the groups
$\GU_n(q)$ for $n \leq 10$, if $\ell > n$. In these cases, most of
the decomposition numbers and the Harish-Chandra series have been
computed by Dudas and Malle~\cite{DuMa}. The Harish-Chandra branching
graphs can be determined from this information using some additional
arguments. The corresponding crystal graphs can be computed with
the GAP3 programs written by one of the authors (see \cite{JaGap3}).

Conjecture~\ref{Uglov} holds for $n = 12$ and $e = 3$ if $\ell \geq 13$.
\subsection{}
\label{MainSeriesLabelling}
There are cases where Conjecture~\ref{Uglov} is known to be true.
\begin{thm*}
Let $0 \leq t < (e - 1)/2$ be an integer, put $r := t(t+1)/2$ and let 
$\lambda := \Delta_t$.

Let $m \in \mathbb{N}$, put $n := r + 2m$ and $G := \GU_n(q)$. Then
$$L := L_{r,m} \cong \GU_r(q) \times \GL_1(q^2)^m$$
is a pure Levi subgroup of~$G$ and~$X_\lambda$ is a cuspidal unipotent 
$kL$-module.

If~$\ell$ is large enough, the unipotent $kG$-module $X_\rho$ lies in
the Harish-Chandra defined by $(L,X_\lambda)$ if and only if
$$\bar{\rho}^{(2)} \in \Phi_{e,m}^{(t +(1-e)/2,0)},$$
where $\Phi_{e,m}^{(t +(1-e)/2,0)}$ denotes the set of Uglov bipartitions 
of~$m$. (See \cite[Definition~$4.4$]{GeJa}; the Uglov bipartitions are simply 
the bipartitions labelling the vertices of the connected component of the 
crystal graph containing $|(-,-), \bc\rangle$.)
\end{thm*}
\begin{prf}
The cuspidal unipotent $KG$-module~$Y_\lambda$ of $\GU_r(q)$ reduces irreducibly 
to the unipotent $kG$-module $X_\lambda$ (see \cite[Theorem~$6.10$]{GeHiMaI}). 
In particular,~$X_\lambda$ is cuspidal. 

Let $\hat{X}_\lambda$ denote the (unique) $\mathcal{O}L$-lattice 
in~$Y_\lambda$. The endomorphism algebra $H_{\mathcal{O}}(L, \hat{X}_\lambda)$ 
is an Iwahori-Hecke algebra over~$\mathcal{O}$ of type~$B_m$ with parameters 
$q^{2t+1}$ and~$q^2$.
By a result of Dipper \cite[Theorem~$4.9$]{DipHom}, the $\ell$-modular 
decomposition matrix of $H_{\mathcal{O}}(L, \hat{X}_\lambda)$ is 
embedded into the decomposition matrix of the unipotent $KG$-modules as 
a submatrix.

By our assumption,~$\ell$ does not divide the order of~$L$ and thus
$X_\lambda$ and $\hat{X}_\lambda$ are projective. It follows that 
$R_L^G( \hat{X}_{\lambda} )$ is projective. The corresponding columns of the
decomposition matrix of~$\mathcal{O}G$ are exactly the columns of the
decomposition matrix of $H_{\mathcal{O}}(L, \hat{X}_\lambda)$. 
Let $\hat{Z}$ be an indecomposable summand of $R_L^G( \hat{X}_\lambda )$ 
and let~$Y_\rho$ be a composition factor of $K \otimes_{\mathcal{O}} \hat{Z}$
with~$\rho$ maximal. Then~$X_\rho$ equals the head of $k \otimes_{\mathcal{O}} 
\hat{Z}$ and thus lies in the Harish-Chandra series defined by 
$(L, X_\lambda )$. Every element of this series arises in this way.

To proceed, we will make use of the notion of a \textit{canonical basic
set} as defined in~\cite[Definition~$3.2.1$]{GeJaHecke}. Applying the 
results of \cite[Section~$3$]{GeckPS}, we obtain the following facts.
Firstly, the Iwahori-Hecke algebra $H_k(L, X_\lambda)$ has a canonical 
basic set with respect to Lusztig's $a$-function on $H_k(L, X_\lambda)$ 
(see \cite[p.~$13$]{GeJaHecke}), if $\ell$ is large enough.
Secondly, this canonical basic set agrees with the canonical basic set
of a suitable specialization of a generic Iwahori-Hecke algebra to an
Iwahori-Hecke algebra $H_K^{(2e)}$ of type~$B_m$, whose parameters are 
powers of a $2e$th root of unity. 
The canonical basic set of $H_k(L, X_\lambda)$ (or rather of the algebra
$H_K^{(2e)}$), is determined in \cite[Theorem~$5.4$, Example~$5.6$]{GeJa}. 
The elements of this canonical basic set are labelled by the set of Uglov 
$m$-bipartitions.

The simple $H_K( L, Y_\lambda )$-modules correspond to the simple constituents 
of $R_L^G( Y_\lambda )$. Arrange the latter by lexicographically decreasing 
labels. By \cite[Theorem~$3.7$]{GeMa} and the results of Lusztig summarized
in \cite[2.2.12]{GeJaHecke}, this ordering corresponds to the ordering of the 
simple modules of $H_K( L, Y_\lambda )$ via Lusztig's $a$-function. Through the 
embedding of the decomposition matrix of $H_{\mathcal{O}}(L, \hat{X}_\lambda)$, 
the members of the canonical basic set thus correspond to the composition 
factors of $R_L^G( Y_\lambda )$ which are at the top of their respective columns 
in the decomposition matrix of~$\mathcal{O}G$. As these top composition factors 
label the $kG$-modules in the Harish-Chandra series of~$kG$ defined by 
$(L,X_\lambda)$, our claim follows.
\end{prf}

\noindent
Theorem~\ref{MainSeriesLabelling} is true without the assumption that 
$t < (e-1)/2$ if Conjecture~\ref{Compatibility} holds. Indeed, in this case 
every unipotent $kG$-module in 
the $(L,X_\lambda)$-series is labelled by a partition with $2$-core~$\Delta_t$. 
Let $\widehat{P(X_\lambda)} \in \mathcal{O}L$-mod denote the projective cover 
of $X_\lambda$. Again by \cite[Theorem~$4.9$]{DipHom}, the decomposition matrix 
of $R_L^G( \widehat{P(X_\lambda)} )$ contains the decomposition matrix of 
$H_{\mathcal{O}}(L, Y_\lambda )$ as a submatrix (with a row of the latter 
labelled by $\mu \vdash_2 m$ corresponding to a row of the former labelled
by $\Phi_t(\mu)$). Let $\hat{Z}$ be an indecomposable summand
of $R_L^G( \widehat{P(X_\lambda)} )$ such that the head of 
$k \otimes_{\mathcal{O}} \hat{Z}$ lies in the Harish-Chandra series defined by 
$(L,X_\lambda)$. Put $Z := K \otimes_{\mathcal{O}} \hat{Z}$, and let~$Y_\rho$ be 
a unipotent composition factor of~$Z$ with~$\rho$ maximal. Then~$X_\rho$ is the 
head of $k \otimes_{\mathcal{O}} \hat{Z}$, and hence $\rho_{(2)} = \Delta_t$.
It follows as in the proof above that $\bar{\rho}^{(2)} \in 
\Phi_{e,m}^{(t +(1-e)/2,0)}$.

\subsection{}
\label{ParameterOf2Core}
Provided Conjecture~\ref{CuspidalityAndECore} is true, we can compute the 
parameters of $H_k(L,X)$ for weakly cuspidal pairs $(L,X)$. We use the 
notation of Theorem~\ref{StructureHLX} in the following.
\begin{prp*}
Suppose that~$X$ lies in a $kM$-block $\mathbf{B}$ whose 
$e$-core equals the $2$-core $\Delta_s$ for some $s \geq 0$. Then 
$Q = q^{2s+1}$.
\end{prp*}
\begin{prf}
By the results summarized in~\ref{CuspidalityAndECore}, the block~$\mathbf{B}$ 
contains a cuspidal simple $KM$-module~$Y$.  By Theorem~\ref{StructureHLX}(a), 
the parameter~$Q$ is equal to the corresponding parameter of the Iwahori-Hecke 
algebra $H_K( L, Y )$. By the results of Lusztig \cite[Section~$5$]{LuszClass}, 
we have $Q = q^{2s+1}$.
\end{prf}

\subsection{}
\label{DifferentGraphs}
If the Conjectures~\ref{GraphComparison} and~\ref{Uglov} are true, 
Proposition~\ref{ParameterOf2Core} implies a compatibility between 
certain connected components of the crystal graph. 

Suppose that $X_\lambda$ is weakly cuspidal, that $\lambda_{(2)} = \Delta_t$,
and that the $e$-core of $\lambda$ equals $\Delta_s$. (The $e$-core of
$\lambda$ should be a $2$-core by Conjecture~\ref{CuspidalityAndECore}.)

Put $r := t(t+1)/2$ and suppose that $n = r + 2m$ and let~$L$ denote the
pure Levi subgroup of~$\GU_n(q)$ isomorphic to $\GU_r(q) \times \GL_1(q^2)^m$.
By Theorem~\ref{StructureHLX} and Proposition~\ref{ParameterOf2Core}, we have
that $H_k(L,X_\lambda$) is an Iwahori-Hecke algebra of type $B_m$ with 
parameters $q^{2s + 1}$ and $q^2$.
According to \cite{GeJa}, the irreducible modules of this Hecke algebra
are labelled by $\Phi_{e}^{(s +(1-e)/2,0)} = 
\cup_{m\geq 0} \Phi_{e,m}^{(s +(1-e)/2,0)}$.
By the generalization
of \cite[Theorem~$2.4$]{GeHiMaII} to weakly cuspidal modules, the elements
of the $(L,X_\lambda)$-Harish-Chandra series of~$kG$ are labelled by these
bipartitions (see also Proposition~\ref{WeakHCSeries}(a)).

On the other hand, by Conjecture~\ref{Uglov}, this Harish-Chandra series
should also be labelled by the set of bipartitions arising from 
$\bar{\lambda}^{(2)}$ by adding a sequence of good nodes with respect
to the charge $(t + (1-e)/2,0)$.

The compatibility of the two labellings is guaranteed by 
Theorem~\ref{ComparingDifferentCrystalGraphs1} below.

\subsection{}
\label{OneFourSeries}
We give an example for the phenomenon discussed above.
Suppose that $e = 3$ and let $L := \GU_4(q) \times \GL_1(q^2)^m$. Then the 
Steinberg $kL$-module
$X_{(1^4)}$ is cuspidal. As the $2$-core of $(1^4)$ is trivial we have $t = 0$.
According to Conjecture~\ref{GraphComparison}, the connected component of the
Harish-Chandra branching graph beginning in $(1^4)$ should coincide, up to some
rank depending on~$\ell$, with the component of the crystal graph corresponding 
to $e = 3$ and charge $(-1,0)$ containing the bipartition $( - ,1^2)$.

The Iwahori-Hecke $H_k(L,X_{(1^4)})$ is of type~$B_m$ with parameters~$q^3$ 
and~$q^2$, as $s = 1$. Its simple modules are labelled by the 
Uglov-bipartitions corresponding to $e = 3$ and charge $(0,0)$.

\subsection{}
\label{WeightOne}
For blocks of $e$-weight~$1$ (for the notions of $e$-core and $e$-weight of
a unipotent $\ell$-block of~$G$ see~\ref{FongSrinivasan}), 
Conjecture~\ref{CuspidalityAndECore} is true.
\begin{thm*}
Let $\mathbf{B}$ be a unipotent $\ell$-block of $\GU_n(q)$ of $e$-weight~$1$.
Then $\mathbf{B}$ contains a weakly cuspidal $kG$-module, if and only if the 
$e$-core of~$\mathbf{B}$ is a $2$-core.
\end{thm*}
\begin{prf}
Suppose first that the $e$-core of~$\mathbf{B}$ is a $2$-core. Then~$\mathbf{B}$
contains a cuspidal simple $KG$-module by the results recalled 
in~\ref{CuspidalityAndECore}. In particular,~$\mathbf{B}$ contains a cuspidal 
unipotent $kG$-module.

Now suppose that the $e$-core of~$\mathbf{B}$ is not a $2$-core. Let 
$s(\mathbf{B})$ denote the Scopes number of~$\mathbf{B}$ (see \cite[7.2]{HiKeI} 
for the definition of $s(\mathbf{B})$). Our assumption implies that 
$s(\mathbf{B}) \geq 1$. Indeed, consider an $e$-abacus diagram (in the sense of
\cite[p.~$78$f]{jake} or \cite[Section~$1$]{fosri3}) for the $e$-core 
of~$\mathbf{B}$. Since the latter is not a $2$-core, there is 
$0 \leq i \leq e - 1$
such that the number of beads on string~$i$ is at least one larger than the 
number of beads on string $i - 2$, if $ 2 \leq i \leq e - 1$, and at least two 
larger than the number of beads on string $ e - 2$ or $e - 1$, if $i = 0$ 
or~$1$, respectively. This exactly means $s(\mathbf{B}) \geq 1$. The Reduction 
Theorem and its consequence \cite[Theorems~$7.10$,~$8.1$]{HiKeI} now imply that
every projective $kG$-module of $\mathbf{B}$ is obtained from Harish-Chandra
induction of a projective $kG$-module of $\GU_{n-2}(q) \times \GL_1( q^2 )$.
In particular, $\mathbf{B}$ contains no weakly cuspidal $kG$-module.
\end{prf}

\subsection{}
\label{EWeightOneCuspidals}
We now determine all partitions $\mu \in \mathcal{P}$ of $e$-weight~$1$ such
that $X_\mu$ is weakly cuspidal. For $0 \leq t \leq (e-1)/2$ let 
$$\mu_{t,e} := (t, t-1, \ldots 3, 2, 1^{e+1}),$$
and for $0 \leq t < (e-1)/2$ let
$$\nu_{t,e} := (t+2, t+1, \ldots , 3, 2, 1^{e -2t - 2}).$$
(we understand $\mu_{0,e} = 1^e$ and $\mu_{1,e} = 1^{e+1}$). For $t = (e-1)/2$, we
also put $\nu_{t,e} := \mu_{t,e}$.
\begin{prp*}
Let $\mu \in \mathcal{P}_n$ have $e$-weight~$1$. Then~$X_\mu$ is weakly cuspidal
if and only if $n = t(t+1)/2 + e$ for some $0 \leq t \leq (e-1)/2$ and $\mu \in 
\{ \mu_{t,e}, \nu_{t,e} \}$.
\end{prp*}
\begin{prf}
Let~$\mathbf{B}$ denote the unipotent $\ell$-block of~$G$ containing~$X_\mu$.

Assume first that~$X_\mu$ is weakly cuspidal. Then, by Theorem~\ref{WeightOne}, 
the $e$-core of~$\mu$ is a $2$-core, $\Delta_t$, say. In particular, $n = 
t(t+1)/2 + e$. As~$\Delta_t$ is an $e$-core, we have $0 \leq t \leq (e - 1)/2$.

By \cite[(6A)]{fosri3}, the partitions $\mu_{t,e}$ and $\nu_{t,e}$ label the unipotent 
$KG$-modules in~$\mathbf{B}$ connected to the exceptional vertex of the Brauer 
tree of~$\mathbf{B}$ (there is only one such if $t = (e-1)/2$).

Assume that $ \mu \not\in 
\{\mu_{t,e}, \nu_{t,e} \}$. Let~$\mu' \in \{\mu_{t,e}, \nu_{t,e} \}$ such that $Y_\mu$
and~$Y_{\mu'}$ lie on the same side of the exceptional vertex in the Brauer
tree of~$\mathbf{B}$. Then~$\mu$ and~$\mu'$ have the same $2$-core~$\Delta_s$,
say, again by \cite[(6A)]{fosri3}. If $\mu' = \mu_{t,e}$, we clearly have $s < t$,
and thus $\Delta_s$ is an $e$-core. If $\mu' = \nu_{t,e}$, then $s = t + 2$, and
$\Delta_s$ is an $e$-core if $e \geq 2t + 5$, and of $e$-weight~$1$ if $e = 2t + 3$. 
In the latter case, $n = t(t+1)/2 + (2t + 3) = (t+2)(t+3)/2 = |\Delta_s|$, and thus
$\mu = \Delta_s = \nu_{t,e}$, a contradiction. Thus in any case~$\Delta_s$ is an 
$e$-core, and so~$X_{\Delta_s}$ is projective. Using \cite[(6A)]{fosri3} once more,
we find that~$X_\mu$ lies in the Harish-Chandra series defined by $(L,X_{\Delta_s})$, 
where~$L$ is the pure standard Levi subgroup of~$G$ corresponding to 
$\GU_{|\Delta_s|}(q)$. In particular,~$X_\mu$ is not weakly cuspidal, contradicting 
our assumption.

Now assume that~$\mu$ is one of~$\mu_{t,e}$ or~$\nu_{t,e}$. Then the $e$-core of~$\mu$
equals~$\Delta_t$, and~$X_\mu$ corresponds to the edge of the Brauer tree 
linking~$Y_\mu$ with the exceptional vertex. By the results summarized 
in~\ref{CuspidalityAndECore}, the exceptional vertex labels cuspidal simple 
$KG$-modules. Thus~$X_\mu$ is cuspidal. This completes our proof.
\end{prf}

More evidence for our conjectures is given in the next section where we prove
some consequences of our conjectures for the crystal graph.

\section{Some properties of the crystal graph}

The conjectures formulated in Section~\ref{Conjectures} imply some combinatorial
properties of the crystal graphs involved. In this final section we prove some 
of these properties. Throughout this section we let~$e$ and~$t$ be non-negative 
integers with~$e$ odd and larger than~$1$. (Contrary to previous usage, the
letter~$k$ no longer denotes a field, but just an integer.)

\subsection{}
\label{SymbolOfVertex}
Following~\cite[6.5.17]{GeJaHecke}, we define a $1$-runner abacus to be a
subset~$\mathfrak{A}$ of~$\mathbb{Z}$ such that $-j \in \mathfrak{A}$ and
$j \not\in \mathfrak{A}$ for all $j \geq n$ and some $0 \neq n \in \mathbb{N}$.
Let~$\mathfrak{A}$ be a $1$-runner abacus. We enumerate the elements
of~$\mathfrak{A}$ by $a_1, a_2, \ldots $ with $a_1 > a_2 > \cdots$. The elements
of $\mathbb{Z} \setminus \mathfrak{A}$ are called \textit{the holes}
of~$\mathfrak{A}$. If we define $\lambda_j$ to be the number of holes
of~$\mathfrak{A}$ less than~$a_j$, $j = 1, 2, \ldots$, then
$\lambda := (\lambda_1, \lambda_2, \ldots )$ is \textit{the partition associated
to}~$\mathfrak{A}$. The \textit{charge} of~$\mathfrak{A}$ is the integer 
$a_1 - \lambda_1$. Let~$n$ be a positive integer such that 
$\{ -j \mid j \geq n \} \subseteq \mathfrak{A}$. Then the number of elements 
of~$\mathfrak{A}$ larger than~$-n$ equals~$n$ plus the charge of~$\mathfrak{A}$.
Moreover, a $\beta$-set for~$\lambda$, in the sense of 
\cite[p.~$2$]{JaYou}, is obtained by 
adding a constant~$d$ to the elements of 
$\mathfrak{A} \setminus \{ -j \mid j \geq n \}$ to make them all non-negative.
Let~$\mathfrak{A}$ 
and~$\mathfrak{A}'$ be $1$-runner abaci with associated partitions~$\lambda$
and~$\lambda'$ and charges~$c$ and~$c'$. Then $\mathfrak{A} = \mathfrak{A}'$
if and only if $\lambda = \lambda'$ and $c = c'$. Also, if $\mathfrak{A}
\subseteq \mathfrak{A}'$ and $|\mathfrak{A}' \setminus \mathfrak{A}| = 1$,
then $c' = c + 1$.

By a \textit{symbol} we mean a pair $\mathfrak{B} := 
(\mathfrak{B}^1,\mathfrak{B}^2)$ of $1$-runner abaci. The components
$\mathfrak{B}^1$ and $\mathfrak{B}^2$ are also called the \textit{first}
and \textit{second row} of $\mathfrak{B}$, respectively. If $\mu^i$ and~$c_i$ 
are the partition associated to~$\mathfrak{B}^i$ and the charge 
of~$\mathfrak{B}^i$, respectively, $i = 1, 2$, we also write $\mathfrak{B} = 
\mathfrak{B}( \mu, \mathbf{c} )$ with $\mu = (\mu^1, \mu^2 )$ and 
$\mathbf{c} = (c_1, c_2)$. Let $\mathbf{c} = (c_1,c_2) \in \mathbb{Z}^2$ and 
let $\mu \in \mathcal{P}^{(2)}$ be a bipartition. Then 
$\mathfrak{B}( \mu, \mathbf{c} )$ can be computed as follows 
(see \cite[2.2]{JL}). Let $\mu = (\mu^1, \mu^2)$ with
$\mu^i = ( \mu^i_{j} )_{j \geq 1}$ and $\mu^i_{j} \geq \mu^i_{j+1} \geq 0$ for 
$j \geq 1$ and $i = 1, 2$. Then $\mathfrak{B}(\mu,\mathbf{c}) = 
( \mathfrak{B}(\mu,\mathbf{c})^1, \mathfrak{B}(\mu,\mathbf{c})^2)$ with 
$\mathfrak{B}(\mu,\mathbf{c})^i := \mathfrak{B}(\mu,\mathbf{c})^i_{j}$, where
$\mathfrak{B}(\mu,\mathbf{c})^i_{j} := \mu^i_{j} - j + c_i + 1$ for $i = 1, 2$ 
and $j \geq 1$. 

\subsection{}
\label{AbacusForPhi}
Put $\bc = (t + (1-e)/2, 0 )$ and let $\mu = (\mu^1, \mu^2)$ be a bipartition. 
To $\mathfrak{B}(\mu,\mathbf{c})$ we associate the $1$-runner abacus 
$$\mathfrak{A}_e(\mu,\mathbf{c}) := 
\{ 2j + e \mid j \in \mathfrak{B}(\mu,\mathbf{c})^1 \}
\cup \{ 2j \mid j \in \mathfrak{B}(\mu,\mathbf{c})^2 \}.
$$
In order to determine the partition associated to 
$\mathfrak{A}_e(\mu,\mathbf{c})$, choose an even positive integer $n = 2m$ such 
that $\{ -j \mid j \geq n - 1 \} \subseteq \mathfrak{A}_e(\mu,\mathbf{c})$ and
put
$$\bar{\mathfrak{A}} := 
\{ x + n \mid x \in \mathfrak{A}_e(\mu,\mathbf{c}), x \geq - n \}.$$
Then~$\bar{\mathfrak{A}}$ is a $\beta$-set for the partition associated to
$\mathfrak{A}_e(\mu,\mathbf{c})$ with $0, 1 \in \bar{\mathfrak{A}}$.
Let 
$$\bar{\mathfrak{A}}^1 := 
\{ (x-1)/2 \mid x \in \bar{\mathfrak{A}}, x \text{\ odd} \}$$
and 
$$\bar{\mathfrak{A}}^2 := 
\{ x/2 \mid x \in \bar{\mathfrak{A}}, x \text{\ even} \}.$$
Then 
$$\bar{\mathfrak{A}}^1 = 
\{ j + (e-1)/2 + m \mid j \in \mathfrak{B}(\mu,\mathbf{c})^1, j \geq -m - (e - 1)/2 \}$$
and 
$$\bar{\mathfrak{A}}^2 = 
\{ j + m \mid j \in \mathfrak{B}(\mu,\mathbf{c})^2, j \geq -m \}.$$
In particular, $\bar{\mathfrak{A}}^i$ is a $\beta$-set for~$\mu^i$, $i = 1, 2$
and $|\bar{\mathfrak{A}}^1| = |\bar{\mathfrak{A}}^2| + t$. The latter equality
follows from the remarks in the first paragraph of~\ref{SymbolOfVertex}.
\begin{lem*}
The partition associated to 
$\mathfrak{A}_e( \mu, \mathbf{c} )$ equals~$\Phi_t(\mu)$.
\end{lem*}
\begin{prf}
Use the notation introduced above. Then 
$|\bar{\mathfrak{A}}| = (|\bar{\mathfrak{A}}^1| + |\bar{\mathfrak{A}}^2|)
\equiv t\,(\text{mod\ } 2)$. Thus $\bar{\mathfrak{A}}$ 
is a $\beta$-set for the partition with $2$-core~$\Delta_t$, and $2$-quotient
(computed with respect to a $\beta$-set with an odd number of elements)
$(\mu^2,\mu^1)$ if~$t$ is odd, and $(\mu^1,\mu^2)$ if~$t$ is even. This
implies our claim.
\end{prf}

\subsection{}
\label{SymbolOfVertex2}
Let $\bc = (t + (1-e)/2,0)$ and let $\mu \in \mathcal{P}^{(2)}$.
We are interested in the operation of deleting $e$-hooks from~$\Phi_t(\mu)$. 
On $\mathfrak{A}_e(\mu,\mathbf{c})$, this amounts to replacing an element 
$y \in \mathfrak{A}_e(\mu,\mathbf{c})$ with $y - e 
\not\in \mathfrak{A}_e(\mu,\mathbf{c})$ by $y - e$.
If~$y$ is odd, this replacement corresponds to the operation of deleting 
$j = (y-e)/2$ from $\mathfrak{B}(\mu,\mathbf{c})^1$ and inserting~$j$ into 
$\mathfrak{B}(\mu,\mathbf{c})^2$. If~$y$ is even, this replacement corresponds 
to the operation of deleting $j = y/2$ from $\mathfrak{B}(\mu,\mathbf{c})^2$ and 
inserting $j - e$ into $\mathfrak{B}(\mu,\mathbf{c})^1$. 
This leads to the following operations on symbols, to which we refer as
\textit{elementary operations}.
\begin{itemize}
\item[(a)] Delete an element~$j$ in the first row, which is not in the
second row, and insert~$j$ in the second row.
\item[(b)] Delete an element~$j$ in the second row, such that~$j - e$
is not in the first row, and insert $j-e$ in the first row.
\end{itemize}
Iterating the two operations we end up with a symbol for which no such operation
is possible. Even though the resulting symbol does not depend on the order in
which we perform these operations, we decide to do the former operation first if 
possible, and always take the largest possible~$j$ so that each step in the 
algorithm is well defined. 
This gives the following elementary operations in a more restrictive sense.
\begin{itemize}
\item[(a$'$)] Delete the largest element~$j$ in the first row, which is not in the 
second row, and insert~$j$ in the second row.
\item[(b$'$)] If every element in the first row is contained in the second row,
delete the largest element~$j$ in the second row, such that~$j - e$ 
is not in the first row, and insert $j-e$ in the first row.
\end{itemize}

\begin{prp*}
Put $\lambda := \Phi_t( \mu )$. Let $\mu' = ( (\mu')^1, (\mu')^2 ) \in 
\mathcal{P}^{(2)}$ and $\mathbf{c}' \in \mathbb{Z}^2$
such that $\mathfrak{B}( \mu',\mathbf{c}' )$ is obtained from
$\mathfrak{B}( \mu,\mathbf{c} )$ by an elementary operation of type~(a) or~(b).

Applying this elementary operation corresponds to removing an $e$-hook from
$\lambda$. Denote by $\lambda'$ the resulting partition, and let $t'$ be
such that $\lambda'_{(2)} = \Delta_{t'}$. Suppose that $\tilde{\mu} =
( \tilde{\mu}^1, \tilde{\mu}^2 )$ is the bipartition such that
$\Phi_{t'}(\tilde{\mu}) = \lambda'$. 

Then $t' = t + 2$, if the elementary operation applied is of type~(b).
If the elementary operation applied is of type~(a), then
$$t' =    \left\{ \begin{array}{ll}
t - 2, & \text{ if } t \geq 2, \\
0, & \text{ if } t = 1, \\
1, & \text{ if } t = 0.
\end{array} \right.$$
Moreover,
                     $$\tilde{\mu} = \left\{ \begin{array}{ll}
                        \mu', & \text{ if $t$ and $t'$ have the same parity,} \\
                        ((\mu')^2,(\mu')^1), & \text{ otherwise,} 
                       \end{array} \right.$$
\end{prp*}
\begin{prf}
Consider a $\beta$-set $\bar{\mathfrak{A}}$ for $\Phi_t(\mu)$ as constructed 
in~\ref{AbacusForPhi}. An elementary operation results in replacing an 
element $x$ of $\bar{\mathfrak{A}}$ by $x - e$ yielding the $\beta$-set 
$\bar{\mathfrak{A}}'$ for~$\lambda'$. (Notice that $\bar{\mathfrak{A}}'$ is
constructed from $\mathfrak{B}( \mu',\mathbf{c}' )$ in the same way as
$\bar{\mathfrak{A}}$ from $\mathfrak{B}( \mu,\mathbf{c} )$.) Moreover,~$x$ is 
even or odd, if the elementary operation is of type~(b) or~(a), respectively. 
In the former case, the number of odd elements of $\bar{\mathfrak{A}}$ increases 
by~$1$, and thus $t' = t + 2$. In the latter case, the number of odd elements of 
$\bar{\mathfrak{A}}$ decreases by~$1$. Hence $t' = t - 2$ if $t \geq 2$, 
$t' = 0$ if $t = 1$, and $t' = 1$ if $t = 0$.

If the parity of~$t$ is the same as that of~$t'$, then the constructions 
of~$\Phi_t(\mu)$ and of $\Phi_{t'}(\tilde{\mu})$ are the same, namely we have 
$\lambda^{(2)}=(\mu^1,\mu^2)$ and $(\lambda')^{(2)} = 
(\tilde{\mu}^1,\tilde{\mu}^2)$ (respectively $\lambda^{(2)}=(\mu^2,\mu^1)$ and 
$(\lambda')^{(2)}= (\tilde{\mu}^2,\tilde{\mu}^1)$) if $t$ is even (respectively 
odd). Therefore, one can read off $\tilde{\mu}$ directly on the symbol 
$\mathfrak{B}(\mu',\bc')$ (or on the $\beta$-sets $\bar{\mathfrak{A}}^i$, 
$i = 1, 2$). It follows that $\tilde{\mu}=\mu'$.

On the contrary, if~$t$ and~$t'$ have different parities (say, without loss of 
generality,~$t$ even and~$t'$ odd), then the construction of 
$\Phi_{t'}(\tilde{\mu})$ requires a permutation, unlike that of $\Phi_t(\mu)$.
Therefore, one needs to permute the components of the bipartition one reads 
off $\mathfrak{B}(\mu',\bc')$, i.e.\ $\tilde{\mu}=((\mu')^2,(\mu')^1)$
\end{prf}

\noindent
As an example, consider the bipartition $\mu = ((5^3,4^2),(6))$, let $e=3$ 
and $t = 5$. Then ${\bc}=(4,0)$ and
$$\mathfrak{B} (\mu,{\bc}) = \left( \begin{array}{llllllll}
\cdots & -2 & - 1 & 6 &   &   &   & \\
\cdots & -2 & - 1 & 4 & 5 & 7 & 8 & 9 
\end{array}
\right).
$$
The associated $1$-runner abacus $\mathfrak{A}_3(\mu,\mathbf{c})$ can
be represented as follows:
$$
\begin{picture}(260,30)(40,-10)
\put( - 9,   5){\tiny $-2$}
\put(   7,   5){\tiny $-1$}
\put(  28,   5){\tiny $0$}
\put(  43,   5){\tiny $1$}
\put(  58,   5){\tiny $2$}
\put(  73,   5){\tiny $3$}
\put(  88,   5){\tiny $4$}
\put( 103,   5){\tiny $5$}
\put( 118,   5){\tiny $6$}
\put( 133,   5){\tiny $7$}
\put( 148,   5){\tiny $8$}
\put( 163,   5){\tiny $9$}
\put( 175,   5){\tiny $10$}
\put( 190,   5){\tiny $11$}
\put( 205,   5){\tiny $12$}
\put( 220,   5){\tiny $13$}
\put( 235,   5){\tiny $14$}
\put( 250,   5){\tiny $15$}
\put( 265,   5){\tiny $16$}
\put( 280,   5){\tiny $17$}
\put( 295,   5){\tiny $18$}
\put( 310,   5){\tiny $19$}
\put( 325,   5){\tiny $20$}
\put( 340,   5){\tiny $21$}
\put(   0,   0){\circle*{5}}                                                     
\put(  15,   0){\circle*{5}}                                                     
\put(  30,   0){\circle*{2}}                                                     
\put(  45,   0){\circle*{5}}                                                     
\put(  60,   0){\circle*{2}}                                                     
\put(  75,   0){\circle*{2}}                                                     
\put(  90,   0){\circle*{2}}                                                     
\put( 105,   0){\circle*{2}}                                                     
\put( 120,   0){\circle*{2}}                                                     
\put( 135,   0){\circle*{2}}                                                     
\put( 150,   0){\circle*{2}}                                                     
\put( 165,   0){\circle*{2}}                                                     
\put( 180,   0){\circle*{2}}
\put( 195,   0){\circle*{5}}
\put( 210,   0){\circle*{5}}
\put( 225,   0){\circle*{5}}
\put( 240,   0){\circle*{2}}                                   
\put( 255,   0){\circle*{2}}                                   
\put( 270,   0){\circle*{2}}                                   
\put( 285,   0){\circle*{5}}                                   
\put( 300,   0){\circle*{2}}                                   
\put( 315,   0){\circle*{5}}                                   
\put( 330,   0){\circle*{2}}                                   
\put( 345,   0){\circle*{5}}                                   
\put( -10, 0.5){\line(1,0){365}}                                     
\end{picture}                                                         
$$
With the notation of~\ref{AbacusForPhi}, taking $n = 2$ we obtain the 
$\beta$-set $\bar{\mathfrak{A}} = \{ 0,1,3,13,14,15,19,21,23\}$ for the 
partition $\lambda := \Phi_5( \mu )$ associated to 
$\mathfrak{A}_3(\mu,\mathbf{c})$. We also have 
$\bar{\mathfrak{A}}^1 = \{ 0, 1, 6, 7, 9, 10, 11 \}$ and 
$\bar{\mathfrak{A}}^2 = \{ 0, 7 \}$, wich are $\beta$-sets for $\mu^1 = 
(5^3,4^2)$ and $\mu^2 = (6)$ respectively. Notice that
$\lambda = (15,14,13,10^3,1)$.
An elementary operation of type~(a$'$) on the symbol yields
$$\mathfrak{B}(\mu',{\bc}') = \left( \begin{array}{llllllll}
\cdots & -2 & - 1 & 6 & 9 &   &   & \\
\cdots & -2 & - 1 & 4 & 5 & 7 & 8 &  
\end{array}
\right),
$$
with $\mu' = ( (5^2,4^2), (8,6) )$ and $\bc' = ( 3, 1 ) )$.
The $1$-runner abacus $\mathfrak{A}_3(\mu',\mathbf{c}')$
cna be pictured as follows:
$$
\begin{picture}(260,30)(40,-10)
\put( - 9,   5){\tiny $-2$}
\put(   7,   5){\tiny $-1$}
\put(  28,   5){\tiny $0$}
\put(  43,   5){\tiny $1$}
\put(  58,   5){\tiny $2$}
\put(  73,   5){\tiny $3$}
\put(  88,   5){\tiny $4$}
\put( 103,   5){\tiny $5$}
\put( 118,   5){\tiny $6$}
\put( 133,   5){\tiny $7$}
\put( 148,   5){\tiny $8$}
\put( 163,   5){\tiny $9$}
\put( 175,   5){\tiny $10$}
\put( 190,   5){\tiny $11$}
\put( 205,   5){\tiny $12$}
\put( 220,   5){\tiny $13$}
\put( 235,   5){\tiny $14$}
\put( 250,   5){\tiny $15$}
\put( 265,   5){\tiny $16$}
\put( 280,   5){\tiny $17$}
\put( 295,   5){\tiny $18$}
\put( 310,   5){\tiny $19$}
\put( 325,   5){\tiny $20$}
\put( 340,   5){\tiny $21$}
\put(   0,   0){\circle*{5}}                                                     
\put(  15,   0){\circle*{5}}                                                     
\put(  30,   0){\circle*{2}}                                                     
\put(  45,   0){\circle*{5}}                                                     
\put(  60,   0){\circle*{2}}                                                     
\put(  75,   0){\circle*{2}}                                                     
\put(  90,   0){\circle*{2}}                                                     
\put( 105,   0){\circle*{2}}                                                     
\put( 120,   0){\circle*{2}}                                                     
\put( 135,   0){\circle*{2}}                                                     
\put( 150,   0){\circle*{2}}                                                     
\put( 165,   0){\circle*{2}}                                                     
\put( 180,   0){\circle*{2}}
\put( 195,   0){\circle*{5}}
\put( 210,   0){\circle*{5}}
\put( 225,   0){\circle*{5}}
\put( 240,   0){\circle*{2}}                                   
\put( 255,   0){\circle*{2}}                                   
\put( 270,   0){\circle*{2}}                                   
\put( 285,   0){\circle*{5}}                                   
\put( 300,   0){\circle*{5}}                                   
\put( 315,   0){\circle*{5}}                                   
\put( 330,   0){\circle*{2}}                                   
\put( 345,   0){\circle*{2}}                                   
\put( -10, 0.5){\line(1,0){365}}                                     
\end{picture}                                                         
$$
We obtain $\bar{\mathfrak{A}}' = \{ 0,1,3,13,14,15,19,20,21\}$, again
using $n = 2$. Next, $(\bar{\mathfrak{A}}')^1 = \{ 0, 1, 6, 7, 9, 10 \}$ 
and $(\bar{\mathfrak{A}})^2 = \{ 0, 7, 10 \}$, wich are $\beta$-sets for 
$(5^2,4^2)$ and $(8,6)$ respectively. The partition associated to 
$\mathfrak{A}_3(\mu',\mathbf{c}')$ is $\lambda' = (13^3,10^3,1)$ which is 
obtained from~$\lambda$ by removing a $3$-hook. We have $\lambda'_{(2)} = 
\Delta_3$, i.e.\ $t' = 3$, and $\Phi_3( \mu' ) = \lambda'$.

\subsection{}
\label{SymbolOfVertex3}
In the following we will make use of the notion of an $e$-period of a symbol
(see \cite[Definition~$2.2$]{JL}) and the concept of totally periodic symbols
(see \cite[Definition~$5.4$]{JL}). Let $|\mu, \bc \rangle$ be a charged
bipartition. In our special situation, an $e$-period of $\mathfrak{B}(\mu,\bc)$
is a sequence $(i_1,k_1), (i_2,k_2), \ldots , (i_e,k_e)$ of pairs of integers
with $2 \geq k_1 \geq k_2 \geq \cdots \geq k_e \geq 1$ such that 
$\mathfrak{B}(\mu,\bc)_{i_l}^{k_l} = m - l + 1$ for some integer~$m$. 
Moreover,~$m$ is the largest element in $\mathfrak{B}(\mu,\bc)^1 \cup 
\mathfrak{B}(\mu,\bc)^2$, and if $m - l + 1 \in \mathfrak{B}(\mu,\bc)^1$ for 
some $1 \leq l \leq e$, then $k_l = 1$. Suppose that $\mathfrak{B}(\mu,\bc)$ 
has an $e$-period $(i_1,k_1), (i_2,k_2), \ldots , (i_e,k_e)$. Then this 
$e$-period is unique and the entries $\mathfrak{B}(\mu,\bc)_{i_l}^{k_l}$ of
$\mathfrak{B}(\mu,\bc)$ are called the elements of the period. Removing these 
elements from $\mathfrak{B}(\mu,\bc)$, we obtain the symbol 
$\mathfrak{B}(\mu',\bc')$ corresponding to a charged bipartition 
$|\mu',\bc'\rangle$ which may or may not have an $e$-period. If iterating this 
procedure ends up in a symbol $\mathfrak{B}(\nu,\bd)$ such that $\nu$ is the empty
bipartition, then $\mathfrak{B}(\mu,\bc)$ is called totally periodic.

By \cite[Theorem~$5.9$]{JL}, the symbol
$\mathfrak{B}(\mu,\bc)$ is totally periodic, if and only if $|\mu,\bc \rangle$
is a highest weight vertex of~$\mathcal{G}_{\bc,e}$. If~$\mathfrak{B}(\mu,\bc)$ 
is totally periodic, then for each entry~$j$ in~$\mathfrak{B}(\mu,\bc)$, there
is a symbol~$\mathfrak{B}'$, obtained from~$\mathfrak{B}(\mu,\bc)$ 
by removing a sequence of $e$-periods, and an $e$-period $(i_1,k_1), \ldots ,
(i_e,k_e)$ of~$\mathfrak{B}'$, such that $j = (\mathfrak{B}')_{i_l}^{k_l}$ for 
some $1 \leq l \leq e$. By a slight abuse of terminology, we say that~$j$ is
contained in the period $(i_1,k_1), \ldots , (i_e,k_e)$ of~$\mathfrak{B}$. 

Let~$\mathfrak{B}'$ denote the symbol obtained from $\mathfrak{B}(\mu,\bc)$
by applying an elementary operation.
\begin{lem*} \label{lemhw}
If $\mathfrak{B}(\mu,\bc)$ is totally periodic, so is ~$\mathfrak{B}'$.
\end{lem*}
\begin{prf}
Suppose first that $\mathfrak{B}'$ is obtained from 
$\mathfrak{B}(\mu,\bc)$ by an elementary operation~(a). Moving~$j$ from row~$1$ 
to row~$2$ transforms the period $(i_1,k_1), \dots, (i_e,k_e)$ containing~$j$ 
into a period $(i'_1,k'_1),\dots,(i'_e,k'_e)$ such that 
$(\mathfrak{B}')_{i'_l}^{k'_l} = \mathfrak{B}(\mu,\bc)_{i_l}^{k_l}$ for 
all~$l$. In particular, $\mathfrak{B}'$ is also totally periodic.

Suppose now that $\mathfrak{B}'$ is obtained from 
$\mathfrak{B}(\mu,\bc)$ by an elementary operation~(b). Deleting~$j$ from 
row~$2$ and inserting $j-e$ in row~$1$ transforms the period $(i_1,k_1), \dots, 
(i_e,k_e)$ containing $j$ into a period $(i'_1,k'_1),\dots,(i'_e,k'_e)$ such 
that $(\mathfrak{B}')_{i'_l}^{k'_l} = 
\mathfrak{B}(\mu,\bc)_{i_l}^{k_l}-1$ for all $l<e-1$ and 
$(\mathfrak{B}')_{i'_e}^{k'_e} = \mathfrak{B}(\mu,\bc)_{i_e}^{k_e}-e$. 
In particular, $\mathfrak{B}'$ is also totally periodic.
\end{prf}

\subsection{}
\label{SteinbergVertex}
Let $G = \GU_n( q )$, and let~$\ell$ and~$e$ be as in~\ref{UnitaryGroups}.
In \cite[Theorem~$8.3$]{GeHiMaI} we have proved that $X_{(1^n)}$ is cuspidal if
and only if~$e$ is odd and divides~$n$ or $n - 1$. This is consistent with
Conjecture~\ref{Uglov}, as will be shown below. Let $\lambda = (1^n)$. Then 
the $2$-core of~$\lambda$ equals $\Delta_t$ with $t = 0$ if~$n$ is even, and 
$t = 1$ if~$n$ is odd. Also $\bar{\lambda}^{(2)} = ( -, 1^m )$ with
$m = \lfloor n/2 \rfloor$; notice that $n = 2m + t$.

\begin{prp*}
Let $e \geq 3$ be an odd integer, let~$m \in \mathbb{N}$ and $t \in \{ 0, 1 \}$.
Put $\bc := (t + (1 - e)/2,0)$. 

Then the vertex $| (-,1^m), \bc \rangle$ of $\mathcal{G}_{\bc,e}$ is a highest 
weight vertex, if and only if $e \mid 2m + t$ or $e \mid 2m  + t - 1$.
\end{prp*}
\begin{prf}
The proof proceeds by induction on~$m$, the case $m = 0$ being clear. 
Assume that $m > 0$ and let $s, s' \in \{ (e - 1)/2 , (e - 3)/2 \}$ with 
$s \neq s'$. The symbol~$\mathfrak{B}$ of
$| (-,1^m), ( -s, 0) \rangle$ equals
$$\mathfrak{B} = \left( \begin{array}{ccccccccc} 
\cdots & -m & 2 - m & 3 - m & \cdots & - s + 1 & - s + 2 & \cdots & 1 \\
\cdots & -m & 1 - m & 2 - m & \cdots & - s  &  & 
\end{array} \right).$$ 
Let~$\mathfrak{B}'$ be the symbol obtained by removing the $e$-period from~$\mathfrak{B}$.
If $m < e - 1$, we find
$$\mathfrak{B}' = \left( \begin{array}{cccccccccc} 
\cdots & -(e - 1) & \cdots & -m & 2 - m  & \cdots & - s\\
\cdots & -(e - 1) 
\end{array} \right),$$
and if $m = e - 1$, we have
$$\mathfrak{B}' = \left( \begin{array}{cccccccccc} 
\cdots & -(e - 1) & - (e-3)  & \cdots & - s\\
\cdots & -(e - 1) 
\end{array} \right).$$
In the latter two cases,~$\mathfrak{B}'$ does not have an $e$-period and thus~$\mathfrak{B}$ is not
totally periodic. On the other hand,~$e$ does not divide one of $2m - 1$, $2m$,
or $2m + 1$, as $1 \leq m \leq e - 1$.

If $m \geq e$, then
$$\mathfrak{B}' = \left( \begin{array}{cccccccccc} 
\cdots & -m  & 2 -m  & 3 -m  & \cdots & - e + 2 & - e + 3 & \cdots - s\\
\cdots & -m  & 1 -m  & 2 -m  & \cdots & - e + 1  &  & 
\end{array} \right).$$
Thus~$\mathfrak{B}'$ is the symbol of $| (-,1^{m -s - 1}), ( -s', 0) \rangle$.
Now~$\mathfrak{B}$ is totally $e$-periodic if and only if~$\mathfrak{B}'$ is totally $e$-periodic.
By induction,~$\mathfrak{B}'$ is totally $e$-periodic if and only if
$e \mid 2m - 2s - 2$ or $e \mid 2m -2s - 3$ in case $s' = (e - 1)/2$, 
respectively if and only if $e \mid 2m - 2s - 2$ or $e \mid 2m -2s - 1$ 
in case $s' = (e-3)/2$.
Suppose first that $s' = (e-1)/2$. Then $s = (e-3)/2$ and thus $2m - 2s - 2 = 
2m + 1 - e$. The claim follows. The other case works analogously.
\end{prf}

\subsection{}
\label{ConsequenceConjectureUglov}
Let $|\mu,\bc\rangle$ be a charged bipartition, put $\mathfrak{B} := 
\mathfrak{B}(\mu.\bc)$ and $\mathfrak{B}^k := \mathfrak{B}(\mu,\bc)^k$ for 
$k = 1, 2 $.

\begin{lem*}
Suppose that $\mathfrak{B}$ is totally $e$-periodic, that $\mathfrak{B}^1 
\subseteq \mathfrak{B}^2$ and that $j - e \in \mathfrak{B}^1$ for all
$j \in \mathfrak{B}^2$ with $j \geq m$ for some $m \in \mathbb{Z}$.

Then for $k = 1, 2$ we have $j - 1 \in \mathfrak{B}^k$ for all $j \in 
\mathfrak{B}^k$ with $j \geq m - e + 1$.
\end{lem*}
\begin{prf}
Let $j \in \mathfrak{B}^k$ with $j - 1 \not\in \mathfrak{B}^k$. Then 
$j - 1 \not\in \mathfrak{B}^1$ and the period of~$\mathfrak{B}$ 
containing~$j$ ends in~$j$. The first element in this period is $j + e - 1$, 
and $j + e - 1 \in \mathfrak{B}^2$. As $j - 1 = j + e - 1 - e$, it follows 
that $j + e - 1 < m$, hence our claim.
\end{prf}

Put $\mathbf{c} = (t + (1-e)/2,0)$. If Conjecture~\ref{Uglov} is true, the 
highest weight vectors of the crystal graph $\mathcal{G}_{\mathbf{c},e}$ label
the weakly cuspidal unipotent $\GU_n(q)$-modules for large enough primes~$\ell$ 
with $e = e(q,\ell)$. More explicitly, a weakly cuspidal $\GU_n(q)$-module 
$X_\lambda$ with $\lambda_{(2)} = \Delta_t$ should be labelled by the highest 
weight vector $|\bar{\lambda}^{(2)}, \mathbf{c} \rangle$. Moreover, 
if~$X_\lambda$ is weakly cuspidal, the $e$-core of~$\lambda$ should be a 
$2$-core by Conjecture~\ref{CuspidalityAndECore}.

Recall that~$\lambda$ with $\lambda_{(2)} = \Delta_t$ and~$\bar{\lambda}^{(2)}$ 
are related by $\lambda = \Phi_t( \bar{\lambda}^{(2)} )$. 

\begin{thm*}
Let the notation be as above. Let $\mu \in \mathcal{P}^{(2)}$ be such that
$|\mu,\mathbf{c}\rangle$ is a highest weight vertex in 
$\mathcal{G}_{\mathbf{c},e}$. Then the $e$-core of $\Phi_t(\mu)$ is a $2$-core.
\end{thm*}
\begin{prf}
Starting with $\mathfrak{B}(\mu,\mathbf{c})$, we apply a sequence of elementary 
operations, until we reach a symbol $\mathfrak{B}'$, which does not allow any 
such operation. Starting with $\mathfrak{A}_e( \mu, \mathbf{c} )$, the 
corresponding sequence of operations results in a $1$-runner abacus 
$\mathfrak{A}'$, such that $y - e \in \mathfrak{A}'$ for all 
$y \in \mathfrak{A}'$. By Lemma~\ref{AbacusForPhi}, the partition associated 
to $\mathfrak{A}'$ is the $e$-core of $\Phi_t(\mu)$.

The symbol~$\mathfrak{B}'$ is totally $e$-periodic by 
Lemma~\ref{SymbolOfVertex3}, and satisfies the assumptions of the above lemma 
for all $m \in \mathbb{Z}$. Hence for $k = 1, 2$, we have $j - 1 \in 
(\mathfrak{B}')^k$ for every $j \in (\mathfrak{B}')^k$. This implies that 
$x - 2 \in \mathfrak{A}'$ for all $x \in \mathfrak{A}'$. In particular, the 
partition associated to $\mathfrak{A}'$ is a $2$-core.
\end{prf}

We now sketch a different proof of the above theorem. Consider, for $s \in 
\mathbb{Z}$, the space of semi-infinite wedge products $\Lambda^{s+\infty/2}$,
as it is defined in \cite[\S$4$]{Uglov}. We do not 
need the precise definition of this space here but we need to know that there 
are three ways to index the elements of its basis (``the semi-infinite ordered 
wedges''):
\begin{itemize}
\item[-] by the set of elements denoted by $|\lambda,s\rangle$ where $\lambda
         \in \mathcal{P}$;
\item[-] by the set of elements denoted by 
        $| \mu, \bc \rangle$, where $\mu \in \mathcal{P}^{(2)}$ and $\bc =
        (c_1,c_2) \in \mathbb{Z}^2$ is such that $c_1 + c_2 = s$. The way to 
        pass from $|\lambda,s\rangle$ to $|\mu,\bc\rangle$ 
        is purely combinatorial;
\item[-]by the set of elements denoted by $|\lambda^{(e)},{\bc}^{(e)} \rangle$ 
        where $\lambda^{(e)}$ is the $e$-quotient of $\lambda$ and ${\bc}^{(e)}
        =(c_1,\ldots,c_e) \in \mathbb{Z}^e$ satisfies $\sum_{i=1}^e = s$ and
        parametrizes the $e$-core of $\lambda$.
   \end{itemize}
Setting $u:=-v^{-1}$, we have three actions of the algebras 
$\mathcal{U}_v'(\widehat{\mathfrak{sl}_e})$, 
$\mathcal{U}_u' (\widehat{\mathfrak{sl}_2})$ and another algebra  
$\mathcal{H}$ (the Heisenberg algebra) on the space $\Lambda^{s+\infty/2}$. 
Moreover these three actions commute and we have the following decomposition 
(see \cite[Theorem 4.8]{Uglov}):
$$\Lambda^{s+\infty/2}=\bigoplus_{{\bc} \in A^2_e (s)}\mathcal{U}_v'
(\widehat{\mathfrak{sl}_e}).\mathcal{H}.\mathcal{U}_u'(\widehat{\mathfrak{sl}_2}).
|(-,-),{\bc} \rangle,$$
where $A^2_e(s)$ is the set of elements ${\bc}\in \mathbb{Z}^2$ such that 
$c_1-c_2\leq e$ and $c_1 + c_2 = s$. In addition, if we fix $\bc$, the 
associated Fock space of level~$2$ is a 
$\mathcal{U}_v'(\widehat{\mathfrak{sl}_e})$-submodule of $\Lambda^{s+\infty/2}$ 
(that is the actions are compatible).

Let $i \in \{ 0, 1 \}$. Denote by $E_i$ and $F_i$ the Chevalley operators of 
$\mathcal{U}_u'(\widehat{\mathfrak{sl}_2})$. Regarding the action of $E_i$ on 
the set of charged bipartitions following Uglov's work, we see that 
$|\mu, {\mathbf{d}} \rangle$ appears in the expansion of 
$E_i.|\lambda,{\bc}\rangle$ if and only if the symbol of 
$|\mu, \mathbf{d} \rangle$ is obtained from the symbol of 
$|\lambda, \bc \rangle$ by one of the two elementary operations~(a) and~(b) 
described in~\ref{SymbolOfVertex2}. This thus gives an algebraic interpretation 
of these transformations on symbols. Moreover, combining this interpretation 
with some properties of the crystal of $\Lambda^{s+\infty/2}$ (see 
\cite[\S~$4.3$]{Uglov}) leads to an alternative proof of the above theorem.

\subsection{}
\label{ComparingDifferentCrystalGraphs}

For a highest weight vertex $|\mu, \bc \rangle$, write $B(\mu,\bc)$ for the 
connected component of $\mathcal{G}_{\bc,e}$ containing $\mu$.
General crystal theory (see \cite{Kashiwara} for instance) ensures that 
$B(\mu,\bc)\simeq B(\nu,\bd)$ as soon as $|\mu,\bc\rangle$ and $|\nu,\bd\rangle$ 
are both highest weight vertices and $\wt(\nu,\bc)=\wt(\mu,\bd)$. Moreover, by 
the characterization (\ref{weight}), the weights of $|\mu, \bc \rangle$ and 
$|\nu, \mathbf{d} \rangle$ coincide if these two charged bipartitions have the 
same reduced $i$-word for all $0 \leq i \leq e - 1$.

From now on, let $|\mu,\bc\rangle$ be a highest weight vertex in 
$\mathcal{G}_{\bc,e}$. Let $|\mu', \bc' \rangle$ be a charged bipartition such
that~$\mathfrak{B}(\mu',\bc')$ is the symbol obtained from 
$\mathfrak{B}(\mu,\bc)$ by applying one of the elementary operations described
in~\ref{SymbolOfVertex2} (a$'$),~(b$'$). By Lemma~\ref{AbacusForPhi},
this implies in particular that 
$\Phi_t(\mu)$ is not an $e$-core.

\begin{lem*}\label{lemiso}
Under the above hypothesis, $|\mu',\bc'\rangle$ is a highest weight vertex and 
there is a crystal isomorphism $B(\mu,\bc) \simeq B(\mu',\bc')$. 
\end{lem*}
\begin{prf}
By Lemma \ref{lemhw}, we know that $\mathfrak{B}(\mu',\bc')$ is totally 
periodic, and thus $|\mu',\bc'\rangle$ is a highest weight vertex by
\cite[Theorem~$5.9$]{JL}. By the discussion at the beginning of this paragraph, 
it remains to show that the reduced $i$-words of $|\mu, \bc \rangle$ and 
$|\mu', \bc' \rangle$ coincide for all $0 \leq i \leq e - 1$. Denote these 
words by $w_i(\mu, \bc)$ and $w_i(\mu', \bc')$. In this proof, we 
use for more clarity the notation $A_k(j)$ (respectively $R_k(j)$) instead of 
simply $A$ (respectively $R$) to encode the addable (respectively removable) 
node of content $j$ lying in component~$k$ of $\mathfrak{B}(\mu,\bc)$.
Note that the contents of the addable and removable nodes of a bipartition are 
the elements $j-1$ and~$j$, respectively, for~$j$ in the corresponding symbol 
(provided~$j$ encodes a non-zero part). In fact, a removable node of content
$j - 1$ corresponds to an element $j \in \mathfrak{B}( \mu, \bc )^k$ such that
$j - 1 \not\in \mathfrak{B}( \mu, \bc )^k$, and an addable node of content~$j$
corresponds to an element $j \in \mathfrak{B}( \mu, \bc )^k$ such that $j + 1
\not\in \mathfrak{B}( \mu, \bc )^k$. Therefore, since an elementary operation 
affects either just one element~$j$ or just~$j$ and $j-e$, the only differences 
that can occur between $w_i(\mu, \bc)$ and $w_i(\mu', \bc')$ are with 
letters~$A$ and~$R$ corresponding to nodes of content $j-1$, $j$, $j-e-1$ and 
$j-e$. We review the only possible changes by enumerating the cases.

Suppose first that we apply the elementary operation~(a$'$), that is to say we 
move~$j$ from row~$1$ of $\mathfrak{B}(\mu,\bc)$ to row~$2$. Moreover,~$j$ is
the largest element in $\mathfrak{B}(\mu,\bc)^1$ for which $j \not\in 
\mathfrak{B}(\mu,\bc)^2$.  Denote by~$l$ the 
largest element of $\mathfrak{B}(\mu,\bc)^2$. To begin with, assume that~$j$ is 
the largest element of $\mathfrak{B}(\mu,\bc)^1$. 

If $j>l$, then~$j$ is the first element of its period, and thus 
$j-1\in\mathfrak{B}(\mu,\bc)^1$. Moreover, either
\begin{itemize}
\item $j>l+1$, in which case the elementary operation takes $A_1(j)$ to $A_2(j)$
and creates an occurence of $R_2(j-1)A_1(j-1)$, which cancels in the 
reduced $i$-word (for $i = j - 1 \text{\ mod\ } e$), or
 \item $j=l+1$, in which case $A_1(j)$ in $\mathfrak{B}(\mu,\bc)$ becomes 
$A_2(j)$ in $\mathfrak{B}(\mu',\bc')$, and $A_2(j-1)$ becomes $A_1(j-1)$.
\end{itemize}

If $j<l$, the following possibilities arise.
\begin{itemize}
\item If $j+1\notin\mathfrak{B}(\mu,\bc)^2$, then again~$j$ is the first element 
of its period, and thus $j-1\in\mathfrak{B}(\mu,\bc)^1$. Moreover, either
\begin{itemize}
\item[$\ast$] $j-1\notin\mathfrak{B}(\mu,\bc)^2$, and $A_1(j)$ becomes $A_2(j)$ 
and $R_2(j-1)A_1(j-1)$ appears, or 
\item[$\ast$] $j-1\in\mathfrak{B}(\mu,\bc)^2$, 
and $A_2(j-1)$ becomes $A_1(j-1)$ and $A_1(j)$ becomes $A_2(j)$.
\end{itemize}
\item If $j+1\in\mathfrak{B}(\mu,\bc)^2$, then either
\begin{itemize}
\item[$\ast$] $j - 1 \in \mathfrak{B}(\mu,\bc)^1$ and 
$j-1\notin\mathfrak{B}(\mu,\bc)^2$, in which case $R_2(j)A_1(j)$ vanishes and 
$R_2(j-1)A_1(j-1)$ appears, or
\item[$\ast$] $j - 1 \in \mathfrak{B}(\mu,\bc)^1$ and 
$j-1\in\mathfrak{B}(\mu,\bc)^2$, in which case $A_2(j-1)$ becomes $A_1(j-1)$ and
$R_2(j)A_1(j)$ vanishes, or
\item[$\ast$] $j - 1 \not\in \mathfrak{B}(\mu,\bc)^1$ and 
$j-1\notin\mathfrak{B}(\mu,\bc)^2$, in which case $R_1(j-1)$ becomes $R_2(j-1)$
and $R_2(j)A_1(j)$ vanishes, or
\item[$\ast$] $j - 1 \not\in \mathfrak{B}(\mu,\bc)^1$ and 
$j-1\in\mathfrak{B}(\mu,\bc)^2$, in which case~$j$ is the last element in its 
period; if $m \geq j + 1$ is the smallest element of $\mathfrak{B}(\mu,\bc)^2$ 
with $m + 1 \not\in \mathfrak{B}(\mu,\bc)^2$, then~$m$ and $j - 1$ are congruent 
modulo~$e$, and $R_1(j-1)A_2(m)$ and $R_2(j)A_1(j)$ vanish.
\end{itemize}
\end{itemize}
Assume now that $j$ is not the largest element of $\mathfrak{B}(\mu,\bc)^1$.
First we consider the case that $j+1\notin\mathfrak{B}(\mu,\bc)^1$. The fact 
that $\mathfrak{B}(\mu,\bc)$ is totally periodic then implies that $j + 1 \in 
\mathfrak{B}(\mu,\bc)^2$ if $j-1\notin\mathfrak{B}(\mu,\bc)^1$ and $j - 1 
\not\in \mathfrak{B}(\mu,\bc)^2$ if $j-1\in\mathfrak{B}(\mu,\bc)^1$ and
$j + 1 \not\in \mathfrak{B}(\mu,\bc)^2$. We obtain the following five subcases.
\begin{itemize}
\item If $j-1\notin\mathfrak{B}(\mu,\bc)^1$, $j+1\in\mathfrak{B}(\mu,\bc)^2$ 
and $j-1\notin\mathfrak{B}(\mu,\bc)^2$, then $R_2(j)A_1(j)$ vanishes and 
$R_1(j-1)$ becomes $R_2(j-1)$.
\item If $j-1\notin\mathfrak{B}(\mu,\bc)^1$, $j+1\in\mathfrak{B}(\mu,\bc)^2$ 
and $j-1\in\mathfrak{B}(\mu,\bc)^2$, then $R_2(j)A_1(j)$ vanishes and 
$R_2(j-1)A_1(j-1)$ appears.
\item If $j-1\in\mathfrak{B}(\mu,\bc)^1$, $j+1\not\in\mathfrak{B}(\mu,\bc)^2$    
and $j-1\notin\mathfrak{B}(\mu,\bc)^2$, then $A_1(j)$ becomes $A_2(j)$ and 
$R_2(j-1)A_1(j-1)$ appears.
\item If $j-1\in\mathfrak{B}(\mu,\bc)^1$, $j+1\in\mathfrak{B}(\mu,\bc)^2$
and $j-1\notin\mathfrak{B}(\mu,\bc)^2$, then $R_2(j)A_1(j)$ vanishes and
$R_2(j-1)A_1(j-1)$ appears.
\item If $j-1\in\mathfrak{B}(\mu,\bc)^k$, $k = 1, 2$ and 
$j+1\in\mathfrak{B}(\mu,\bc)^2$, then $R_2(j)A_1(j)$ vanishes and $A_2(j-1)$ 
becomes $A_1(j-1)$.
\end{itemize}
If $j+1\in\mathfrak{B}(\mu,\bc)^1$, then $j+1\in\mathfrak{B}(\mu,\bc)^2$ 
otherwise $j$ would not be moved.
\begin{itemize}
\item If $j-1\notin\mathfrak{B}(\mu,\bc)^2$ and 
$j-1\notin\mathfrak{B}(\mu,\bc)^1$, then $R_1(j-1)$ becomes $R_2(j-1)$ and 
$R_2(j)$ becomes $R_1(j)$.
\item The case $j-1\in\mathfrak{B}(\mu,\bc)^2$ and 
$j-1\notin\mathfrak{B}(\mu,\bc)^1$ can not occur as $\mathfrak{B}(\mu,\bc)$ is
totally periodic.
\item If $j-1\notin\mathfrak{B}(\mu,\bc)^2$ and $j-1\in\mathfrak{B}(\mu,\bc)^1$, 
then $R_2(j)$ becomes $R_1(j)$ and $R_2(j-1)A_1(j-1)$ appears.
\item If $j-1\in\mathfrak{B}(\mu,\bc)^2$ and $j-1\in\mathfrak{B}(\mu,\bc)^1$, 
then $R_2(j)$ becomes $R_1(j)$ and $A_2(j-1)$ becomes $A_1(j-1)$.
\end{itemize}

Suppose now that we apply operation~(b$'$), that is to say, that we delete~$j$ 
from $\mathfrak{B}(\mu,\bc)^2$ and insert $j-e$ in $\mathfrak{B}(\mu,\bc)^1$. 
This implies in particular that all elements of $\mathfrak{B}(\mu,\bc)^1$ are in 
$\mathfrak{B}(\mu,\bc)^2$. Again, assume first that~$j$ is the largest element 
of $\mathfrak{B}(\mu,\bc)^2$. As $\mathfrak{B}(\mu,\bc)$ is totally 
$e$-periodic, $j - 1$ appears in $\mathfrak{B}(\mu,\bc)$, hence 
$j - 1 \in \mathfrak{B}(\mu,\bc)^2$. Denote by~$l$ the largest element of  
$\mathfrak{B}(\mu,\bc)^1$.
Suppose first that $j-e>l$.
\begin{itemize}
\item If $j-e>l+1$, then $A_2(j)$ becomes $A_1(j-e)$, and $R_1(j-1-e)A_2(j-1)$ 
appears.
\item If $j-e=l+1$, then  $A_2(j)$ becomes $A_1(j-e)$ and  $A_1(j-e-1)$ becomes 
$A_2(j-1)$.
\end{itemize}
Now assume that $j-e<l$. Note that in this case $j - e + 1 \in 
\mathfrak{B}(\mu,\bc)^1$. ndeed,~$j$ is the first element in the period of
$\mathfrak{B}(\mu,\bc)$, and $j - e + 1$ the last. As $l \geq j - e + 1$
and~$l$ lies in the first row, so does~$j - e + 1$.
\begin{itemize}
\item If $j-e-1\notin\mathfrak{B}(\mu,\bc)^1$, then $R_1(j-e)A_2(j)$ vanishes 
and $R_1(j-e-1)A_2(j-1)$ appears.
\item If $j-e-1\in\mathfrak{B}(\mu,\bc)^1$, then $A_1(j-e-1)$ becomes $A_2(j-1)$ 
and $R_1(j-e)A_2(j)$ vanishes.
\end{itemize}
Finally, assume that~$j$ is not the largest element of $\mathfrak{B}(\mu,\bc)^2$
and let~$j'$ denote the largest element of $\mathfrak{B}(\mu,\bc)^2$. Then 
$j' - e \in \mathfrak{B}(\mu,\bc)^1$, as our operation of type (b$'$) always
moves the largest possible element. Hence $l \geq j' - e > j-e$. Now~$j$ is the 
largest element of $\mathfrak{B}(\mu,\bc)^2$ such that $j - e$ is not in 
$\mathfrak{B}(\mu,\bc)^1$. By Lemma~\ref{ConsequenceConjectureUglov},
this implies that for $k = 1, 2$ and
every $r > j - e + 1$ we have $r - 1 \in \mathfrak{B}(\mu,\bc)^k$ if $r \in 
\mathfrak{B}(\mu,\bc)^k$. Hence all integers in the interval $[j - e + 1, j']$
and $[j - e + 1, l]$ are contained in $\mathfrak{B}(\mu,\bc)^2$ and 
$\mathfrak{B}(\mu,\bc)^1$, respectively. This implies in particular that
$j - e \in \mathfrak{B}(\mu,\bc)^2$ as otherwise the element $j - e + 1$ of the 
second row must be the last element in its period. But then the element 
$j - e + 1$ of the first row must lie in an earlier period, which 
is impossible. This leaves to check the following possibilities.
\begin{itemize}
\item If $j - e - 1 \in \mathfrak{B}(\mu,\bc)^1$ it is also contained in
$\mathfrak{B}(\mu,\bc)^2$, and $R_1( j - e )$ becomes $R_2( j )$ and $A_1(j-e-1)$ 
becomes $A_2( j - 1 )$.
\item If $j - e - 1 \in \mathfrak{B}(\mu,\bc)^2$ and $j - e - 1 \not\in 
\mathfrak{B}(\mu,\bc)^1$, then $R_1( j - e )$ becomes $R_2( j )$ and 
$R_1(j-e-1)A_1(j-1)$ appear.
\item If $j - e - 1 \not\in \mathfrak{B}(\mu,\bc)^2$ and $j - e - 1 \not\in 
\mathfrak{B}(\mu,\bc)^1$, then $R_1( j - e )$ becomes $R_2( j )$ and 
$R_1(j-e-1)A_1(j-1)$ appear.
\end{itemize}

In each case, we see that $w_i(\mu) = w_i(\mu')$, for all $i=1,\dots,e-1$.
\end{prf}

\noindent
We record a first consequence of the above lemma. Let~$t' \in \mathbb{N}$ 
and $\tilde{\mu} \in \mathcal{P}^{(2)}$ be such that $\lambda' := 
\Phi_{t'}(\tilde{\mu})$ equals the partition obtained from $\lambda := 
\Phi_t(\mu)$ be removing the $e$-hook which corresponds to the elementary 
operation transforming $\mathfrak{B}(\mu,\bc)$ into $\mathfrak{B}(\mu',\bc')$. 
(See Proposition~\ref{SymbolOfVertex2} how to compute~$t'$ and~$\tilde{\mu}$.)
Suppose that~$t$ and~$t'$ have the same parity and put $\tilde{\bc} := 
(t' + (1-e)/2,0)$. Then $\tilde{\mu} = \mu'$ and $\tilde{\bc}$ is obtained
from $\bc'$ by adding or subtracting~$1$ to each of its components. By 
definition of the crystal graph, it is clear that translating each component
of the charge by some fixed integer, results in the same graph with an overall 
translation of the labels of the arrows. In particular, 
$| \tilde{\mu}, \tilde{\bc} \rangle$ is a highest weight vertex.
\begin{cor*}
Suppose that Conjecture~\ref{GraphComparison} is true. Then the Harish-Chandra
branching graphs corresponding to the weakly cuspidal modules $X_\lambda$
and $X_{\lambda'}$ are isomorphic (up to some rank).
\end{cor*}
\begin{prf}
It follows from the considerations preceeding the corollary, that 
$B( \mu, \bc )$ and $B( \tilde{\mu}, \tilde{\bc} 0$ 
are isomorphic up to a global shift of the arrow labels.
\end{prf}

\noindent
This corollary shows that the validity of Conjecture~\ref{GraphComparison} 
would yield a remarkable connection between the Harish-Chandra theory of 
unitary groups of odd and even degrees.

\subsection{}
\label{ComparingDifferentCrystalGraphs1}
We finally prove a property of the crystal graph which is implied by the
considerations in~\ref{DifferentGraphs}. Let $\lambda \in \mathcal{P}$ with
$\lambda_{(2)} = \Delta_t$ and $\bar{\lambda}^{(2)} = \mu$. Put $\bc := 
(t + (1-e)/2,0)$. Assume that $|\mu, \mathbf{c}\rangle$ is a highest weight 
vector in $\mathcal{G}_{\mathbf{c},e}$. By 
Theorem~\ref{ConsequenceConjectureUglov}, the $e$-core of~$\lambda$ is a 
$2$-core, $\Delta_s$, say, for some non-negative integer~$s$. Put 
$\mathbf{s} := (s + (1-e)/2,0)$.

\begin{thm*}
With the notation introduced above, there is a graph isomorphism
$$B(\mu,\bc) \simeq B((-,-),\bs),$$
up to a shift of the labels of the arrows.
\end{thm*}
\begin{prf}
We apply the algorithm used to compute the $e$-core of $\lambda = \Phi_t(\mu)$ 
described in the proof of Theorem~\ref{ConsequenceConjectureUglov}. Applying a 
sequence of elementary operations of types~(a$'$) and~(b$'$) to 
$\mathfrak{B}( \mu, \bc )$, we end up with the symbol 
$\mathfrak{B}( (-,-), \bd )$ for some charge $\bd=(d_1,d_2)$. 

We may as well apply the corresponding sequence of moves to the $\beta$-set 
$\bar{\mathfrak{A}}$ for $\lambda = \Phi_t(\mu)$ as constructed 
in~\ref{AbacusForPhi}. This results in a $\beta$-set $\bar{\mathfrak{A}}'$ 
for~$\Delta_s$. The number of odd elements of~$\bar{\mathfrak{A}}$
exceeds its number of even elements by $t = c_1 - (1-e)/2 - c_2$. If the number
of odd elements of~$\bar{\mathfrak{A}}'$ is not smaller than the number of its
even elements, the difference between the two numbers equals~$s$. Otherwise,
there are $s + 1$ more even numbers in~$\bar{\mathfrak{A}}'$ than odd ones.
An operation of type~(a$'$) decreases the first component of the current charge 
by~$1$ and increases the second component by~$1$. The corresponding move on the 
$\beta$-set replaces an odd number by an even one. The analogous remarks apply
for elementary operations of type~(b$'$). We thus find
\begin{equation} \label{charge1} d_1-d_2  = s + (1-e)/2\end{equation}
or
\begin{equation}\label{charge} d_2-d_1  = s + (1+e)/2.\end{equation}

By Lemma~\ref{ComparingDifferentCrystalGraphs}, we have a crystal isomorphism 
$$B(\mu,\bc) \simeq B((-,-),\bd).$$
If we set $\bd'=(d_2-e,d_1)$, we also have a crystal isomorphism
$$B(\mu,\bc) \simeq B((-,-),\bd')$$
(see \cite[6.2.9, 6.2.17]{GeJaHecke}).
By the remark preceding Corollary~\ref{ComparingDifferentCrystalGraphs}, we obtain
$$B(\mu,\bc) \simeq B((-,-),(d_1-d_2,0))$$
and 
$$B(\mu,\bc) \simeq B((-,-),(d_2-e-d_1,0))$$
up to an overall shift of the labels of the arrows. Applying 
Identities~(\ref{charge1}) respectively~(\ref{charge}), we see that 
$s + (1-e)/2$ equals $d_1 - d_2$ in the first case and $d_2 - e - d_1$ in the 
second. This concludes our proof.
\end{prf}

Note that there should be a way to relate these elementary crystal isomorphisms
with the so-called canonical crystal isomorphism of~\cite{Gerber}.

\subsection{}
\label{DifferentBlocks}
Put $\bc := (t + (1-e)/2,0).$ Let $\mu = (\mu^1,\mu^2)$ be a bipartition. For 
$0 \leq j \leq e - 1$, let $\widetilde{f}_{j}$ denote the associated Kashiwara 
operator on $\mathcal{G}_{\bc,e}$ (see~\ref{Fockspace}).

\begin{prp*}
Let $0 \leq j_1 \neq j_2 \leq e - 1$. Suppose that
$\widetilde{f}_{j_i} . |\mu, \bc \rangle \neq 0$ for $i = 1, 2$.  Write 
$\widetilde{f}_{j_i} . |\mu, \bc \rangle = | \nu_i, \bc \rangle$, $i = 1, 2$. 
Then the $e$-cores of $\Phi_t( \nu_1 )$ and of $\Phi_t( \nu_2 )$ are distinct.
\end{prp*}
\begin{prf}
Let $0 \leq j \leq e - 1$. First note that if 
$\widetilde{f}_{j} . |\mu, \bc \rangle \neq 0$ then
\begin{enumerate}
\item $\mathfrak{B}(\mu,\mathbf{c})^{1}=\mathfrak{B}(\widetilde{f}_{j} .
\mu,\mathbf{c})^{1}$ and $\mathfrak{B}(\widetilde{f}_{j} . \mu,\mathbf{c})^{2}=
\mathfrak{B}(\mu,\mathbf{c})^{2} \cup \{k\} \setminus \{k-1\}$
for $k \in \mathbb{Z}$ such that $k \equiv j\,(\textrm{mod\ }e)$, or
\item $\mathfrak{B}(\mu,\mathbf{c})^{2}=\mathfrak{B}(\widetilde{f}_{j} .
\mu,\mathbf{c})^{2}$
and $\mathfrak{B}(\widetilde{f}_{j} . \mu,\mathbf{c})^{1}=
\mathfrak{B}(\mu,\mathbf{c})^{1}\cup \{k\}\setminus \{k-1\}$
for $k\in \mathbb{Z}$ such that $k\equiv j\,(\textrm{mod\ }e)$.
\end{enumerate}
We have seen in~\ref{ConsequenceConjectureUglov} how to compute the $e$-cores of
$\Phi_t (\nu_i )$, $i = 1, 2$. In this procedure, some of the elements $x$ in
$\mathfrak{B}(\nu_i,\mathbf{c})$, $i = 1, 2$,
must be replaced by $x-k.e$ for some $k\in \mathbb{N}$.
If the $e$-core of $\Phi_t (\nu_1 )$ equals the $e$-core of $\Phi_t (\nu_2 )$,
this implies that at the end of these procedures, we obtain the same symbols. 
However, this is impossible as $j_1 \not\equiv j_2\,(\textrm{mod\ }e)$.
\end{prf}

\begin{cor*}
Suppose that Conjecture~\ref{GraphComparison} is true. Let~$X$ be a unipotent 
$k\GU_n(q)$-module. Then, if $\ell$ is large enough, any two non-isomorphic
simple submodules of $R_{\GU_n(q)}^{\GU_{n+2}(q)}( X )$ lie in distinct
$\ell$-blocks.
\end{cor*}
\begin{prf}
By Conjecture~\ref{GraphComparison}, the non-isomorphic simple submodules
of $R_{\GU_n(q)}^{\GU_{n+2}(q)}( X )$ correspond to two distinct directed
edges in a suitable crystal graph. By the Proposition, the corresponding
partitions have distinct $e$-cores, and thus the unipotent modules
labelled by these partitions are in distinct $\ell$-blocks.
\end{prf}

\section*{Acknowledgements}

Above all, we thank Meinolf Geck and Gunter Malle for inspiring this work 
in many ways and for their constant interest. We also thank Klaus Lux and 
Pham Huu Tiep for inviting the second author to the University of Arizona. 
During his visit parts of this paper were written. We thank G\"otz Pfeiffer 
for helpful discussions regarding Lemma~\ref{PureLevis}. Finally, we
very much thank the Isaac Newton Institute for Mathematical Sciences for
inviting two of us to the special programme Algebraic Lie Theory (12 January 
to 26 June 2009), where the ground for the conjectures in this paper was lead.

\begin{landscape}
\begin{table}[f]
\caption{\label{BranchingGraphExample}The Harish-Chandra branching graph for $\GU_n(q)$,
$1 \leq n \leq 7$ odd, $\ell \mid q^2 -q + 1$, $\ell > 7$}
\hrule
\smallskip
$
\begin{xy}
\xymatrix{
& & 1 \ar[d] & & & & & \\
& & 3 \ar[d]\ar[dl]\ar[dr] & & & 1^3 \ar[dl]\ar[dr] & & \\
& 5\ar[dl]\ar[d] & 31^2 \ar[dl]\ar[d] & 32 \ar[dl]\ar[d] & 2^21 \ar[dr]\ar[d] & & 1^5 \ar[d]\ar[dr] & \\
 7 & 51^2 & 32^2 & 52 & 3^21 & 421 & 31^4 & 2^21^3 & 1^7 & 321^2
}
\end{xy}
$
\hrule
\smallskip
\end{table}

\vfill

\bigskip
\bigskip
\begin{table}[f]
\caption{\label{BranchingGraphExample2}The Harish-Chandra branching graph for $\GU_n(q)$,
$1 \leq n \leq 7$ odd, $\ell \mid q^2 -q + 1$, $\ell > 7$, cont.}
\hrule
\smallskip

$
\begin{xy}
\xymatrix{
 & & 21 \ar[dl]\ar[dr] & & &  \\
 &  41 \ar[dl]\ar[d] & &  21^3 \ar[d]\ar[dr] & &  \\
  61 & 43 & & 41^3 & 21^5 & 2^31
}
\end{xy}
$
\hrule
\smallskip
\end{table}

\end{landscape}

\begin{landscape}
\begin{table}[f]
\caption{\label{CrystalGraphExample}The crystal graph $\mathcal{G}_{(0,0),3}^{\leq 3}$}
\hrule
\smallskip
$
\begin{xy}
\xymatrix{
& & -.- \ar[d] & & & & & \\
& & 1.- \ar[d]\ar[dl]\ar[dr] & & & -.1 \ar[dl]\ar[dr] & & \\
& 2.- \ar[dl]\ar[d] & 1.1 \ar[dl]\ar[d] & 1^2.- \ar[dl]\ar[d] & -.2 \ar[dr]\ar[d] & & -.1^2 \ar[d]\ar[dr] & \\
 3.- & 2.1 & 1^2.1 & 21.- & 1.2 & -.3 & 1.1^2 & -.21
& -.1^3 & 1^3.-
}
\end{xy}
$
\hrule
\smallskip
\end{table}
\end{landscape}

\end{document}